\documentclass[12pt]{article}

\usepackage{epsfig,epsfig}
\usepackage{epstopdf}
\usepackage{amsmath}
\usepackage{mathrsfs}
\usepackage{amsfonts}
\usepackage{bm}
\usepackage{amssymb}
\usepackage{amsthm}
\usepackage{stmaryrd}
\usepackage{booktabs}
\usepackage{color}
\usepackage{multirow}
\usepackage{graphicx}
\usepackage{subfigure}
\usepackage{float}
\usepackage{appendix}
\usepackage{tikz}
\usetikzlibrary{calc}
\usetikzlibrary{matrix}
\usepackage{mathtools}
\usepackage{algorithm,algpseudocode,float}
\usepackage{lipsum}


\usepackage{lineno}
\usepackage{srcltx,graphicx}

\newtheorem{thm}{Theorem}[section]

\newtheorem{rem}[thm]{Remark}
\newtheorem{example}[thm]{Example}

\usepackage{hyperref}

\numberwithin{equation}{section} \topmargin=-3.0cm \oddsidemargin=1cm
\begin{document}
\baselineskip=1.5pc

\title{\textbf{A fifth-order finite difference HWENO scheme combined with limiter for hyperbolic conservation laws}}
\author{Min Zhang\footnote{School of Mathematical Sciences, Peking University, Beijing 100871, China. E-mail: minzhang@math.pku.edu.cn. The research of this author is partially supported by the Postdoctoral Science Foundation of China (Grant Nos. 2022M710229).}
~and~
Zhuang Zhao\footnote{School of Mathematical Sciences and Institute of Natural Sciences, Shanghai Jiao Tong University, Shanghai 200240, China. E-mail: zzhao-m@sjtu.edu.cn. The research of this author is partially supported by the Postdoctoral Science Foundation of China (Grant Nos. 2021M702145).}
}

\date{}
\maketitle
\begin{abstract}
In this paper, a simple fifth-order finite difference Hermite WENO (HWENO) scheme combined with limiter is proposed for one- and two- dimensional hyperbolic conservation laws. The fluxes in the governing equation are approximated by the nonlinear HWENO reconstruction which is the combination of a quintic polynomial with two quadratic polynomials, where the linear weights can be artificial positive numbers only if the sum equals one. And other fluxes in the derivative equations are approximated by high-degree polynomials directly. For the purpose of controlling spurious oscillations, an HWENO limiter is applied to modify the derivatives. Instead of using the modified derivatives both in fluxes reconstruction and time discretization as in the modified HWENO scheme (J. Sci. Comput., 85:29, 2020), we only apply the modified derivatives in time discretization while remaining the original derivatives in fluxes reconstruction.
Comparing with the modified HWENO scheme, the proposed HWENO scheme is simpler, more accurate, efficient and higher resolution. In addition, the HWENO scheme has a more compact spatial reconstructed stencil and greater efficiency than the classical fifth-order finite difference WENO scheme of Jiang and Shu. Various benchmark numerical examples are presented to show the fifth-order accuracy, great efficiency, high resolution and robustness of the proposed HWENO scheme.
\end{abstract}

\noindent\textbf{The 2010 Mathematics Subject Classification:} 65M60, 35L65

\vspace{5pt}


\normalsize \vskip 0.2in
\newpage

\section{Introduction}
\label{sec:introduction}
In this paper, we develop a genuine fifth-order finite difference Hermite weighted essentially non-oscillatory (HWENO) scheme for one- and two- dimensional hyperbolic conservation laws. 
Instead of only using the information of the solution in the weighted essentially non-oscillatory (WENO) scheme, the HWENO scheme uses both the information of the solution and its first-order derivatives/moments.
The HWENO scheme gains a more compact stencil in the spatial reconstruction than the WENO scheme on the same order accuracy.
It is well known that the WENO scheme is a particularly powerful numerical tool for the simulation of hyperbolic conservation laws, which was first constructed by Liu, Osher and Chan on the basis of essentially non-oscillatory (ENO) schemes \cite{heoc1,ho1} in 1994.  Since then, Jiang and Shu developed a fifth-order finite difference WENO (WENO-JS) scheme \cite{js} in 1996.
In \cite{js}, the authors gave a general framework for the definition of smoothness indicators and nonlinear weights which is widely used in the subsequent advanced WENO schemes, e.g., \cite{bgsAO,ccd,cd1,hs,lpr,ZStetra,ZQd,ZSWENOMR}, and more detailed reviews for ENO and WENO schemes can refer to \cite{s3,shu}.

The WENO scheme uses the information of the solution in the target cell and its neighbor cells to obtain high-order accuracy, therefore, a higher order WENO scheme will lead to the stencil wider. To resolve this problem, Qiu and Shu \cite{QSHw1} developed a one-dimensional fifth-order finite volume HWENO scheme based on both the information of the solution and its first-order derivative, which only needs the immediate neighbor values in the spatial reconstruction. Meanwhile, for stability, the different reconstructed polynomials are constructed to discretize the fluxes in the original governing and the derivative equations, respectively.
However, this method is not enough to maintain stability and robustness, such as it obtains poor results for the double Mach and the step forward problems in the later two-dimensional work \cite{QSHw}.
Later, Capdeville \cite{Capdeville} developed a finite volume Hermite central WENO scheme. Liu and Qiu \cite{LiuQ1} developed a fifth-order finite difference HWENO scheme in one dimension, unfortunately, it only has the fourth-order accuracy in two dimensions due to the mixed derivatives.
Ma and Wu \cite{MW} developed a compact HWENO scheme by solving the derivatives using the compact difference method.
Recently, Zhao et al. \cite{ZhaoZhuang-2020JSC-FD} developed a genuine fifth-order modified finite difference HWENO (M-HWENO) scheme in one and two dimensions. In \cite{ZhaoZhuang-2020JSC-FD}, the authors modified the derivatives of the solution by a high-order Hermite limiter to control the derivatives near discontinuities and improve the stability of the scheme, while used one set of stencils in the reconstruction, which is different from \cite{LiuQ1,QSHw1,QSHw}. Li et al. \cite{LiMRHW1} developed a multi-resolution HWENO scheme with unequal stencils, but the scheme only has the fourth-order accuracy in two dimensions. For more HWENO schemes, the interested reader can refer to \cite{DBTM,MW,TLQ,weHTENO,ZA,ZQ20,ZQHW} and the references therein.

In this paper, we develop a simple fifth-order finite difference HWENO scheme combined with limiter (denoted as L-HWENO) for one- and two- dimensional hyperbolic conservation laws following the idea of the M-HWENO scheme \cite{ZhaoZhuang-2020JSC-FD}. Instead of using the modified derivatives both in fluxes reconstruction and time discretization as in \cite{ZhaoZhuang-2020JSC-FD}, we only apply the modified derivatives in time discretization while remaining the original derivatives in fluxes reconstruction.
The fluxes in the governing equation are approximated by the nonlinear HWENO reconstruction which is the combination of a quintic polynomial with two quadratic polynomials, where the linear weights can be artificial positive numbers only if the sum equals one. And  other fluxes are approximated by high-degree polynomials directly, which leads to the fact that the reconstruction of the fluxes for derivative equations is linear.  To improve the robustness/stability of the proposed HWENO scheme, the derivatives are modified by a fifth-order HWENO limiter, which is the combination of a quartic polynomial with two linear polynomials using the same technique  in the reconstruction.

Comparing with the M-HWENO scheme \cite{ZhaoZhuang-2020JSC-FD}, the proposed L-HWENO scheme has three main advantages: one is that the modification for derivatives only acts on the time discretization, which makes there no need to storage the original and modified derivatives in the computation as in \cite{ZhaoZhuang-2020JSC-FD}.
The second one is that all the fluxes in the derivative equations can be approximated by high-degree reconstructed polynomials directly, which significantly simplify the algorithm and improve the computational efficiency.
The last one is that both the HWENO limiter and the HWENO spatial reconstruction are based on the combination of a high-degree polynomial with two lower-degree polynomials convexly, which makes the linear weights be any artificial positive constants (their sum equal to one).
It is worth pointing out that the limiter plays an important role to improve stability and keep high resolution, whereas lacking this procedure would lead to instability  in two dimensions  even for a linear problem (cf. Example \ref{Euler-2d}). Meanwhile, different linear weights in the limiter would impact the resolution near discontinuities, while the linear weights in the spatial reconstruction have a slight effect (cf. Example \ref{blast-1d}).

For the spatial reconstruction, the proposed L-HWENO scheme uses a more compact stencil than the same order finite difference WENO-JS scheme \cite{js}. To be specific, the L-HWENO scheme only needs a compact three-point stencil while the WENO-JS scheme needs a five-point stencil in the reconstructions for achieving fifth-order accuracy.
Although the L-HWENO scheme needs to solve the derivative equations which adds extra computational costs into the algorithm, the L-HWENO scheme is more efficient than the WENO-JS scheme in the sense that the former leads to a smaller error than the latter for a fixed amount of the CPU time (cf. Section \ref{sec:numerical}). Note that the efficiency of the M-HWENO scheme \cite{ZhaoZhuang-2020JSC-FD} and the WENO-JS scheme \cite{js} is neck and neck.

The organization of the paper is as follows. In Section \ref{sec:HWENO}, the detailed implementation algorithm of the HWENO scheme combined with limiter is presented in one and two dimensions. In Section \ref{sec:numerical}, various benchmark numerical examples are tested to show the numerical accuracy, great efficiency, high resolution and robustness of the proposed scheme. Concluding remarks are given in Section \ref{sec:conclusions}.

\section{Fifth-order finite difference L-HWENO scheme}
\label{sec:HWENO}

In this section, we present a simple fifth-order finite difference HWENO scheme with limiter (L-HWENO), which combines a high-degree polynomial with two lower-degree polynomials convexly in the spatial reconstruction and limiter, where the associated linear weights both can be chosen as artificial positive number with their sum equals one. Note that the fluxes in the governing equation are approximated by nonlinear HWENO reconstructions, while other fluxes in the derivative equations are approximated by high-degree polynomials directly.

\subsection{Fifth-order finite difference HWENO scheme}
\label{sec:HWENO-1d-u}

For the simplicity of algorithm description, we focus on the scalar equation. The extension to a system of equations is straightforward.
We first consider a one-dimensional scalar hyperbolic conservation laws
\begin{equation}\label{EQ}
\begin{cases}
u_t + f(u)_x=0, \\
u(x,0)=u_0(x).
\end{cases}
\end{equation}
In the finite difference framework, the computing domain is divided by uniform meshes $I_{i} =[x_{i-\frac12},x_{i+\frac12}]$, $\Delta x=x_{i+\frac12}- x_{i-\frac12}$
and $x_i=\frac12( x_{i-\frac12}+x_{i+\frac12})$ is the center of $I_i$. To design an HWENO scheme, we first add the derivative equation of  (\ref{EQ}), having
\begin{equation}\label{EQ1}
\begin{cases}
u_t + f(u)_x=0, \quad u(x,0)=u_0(x),\\
v_t + h(u,v)_x=0, \quad v(x,0)=v_0(x),
\end{cases}
\end{equation}
where $v = u_x$, $h(u,v)=f(u)_x=f'(u)u_x=f'(u)v$. Then, its semi-discrete finite difference HWENO scheme is
\begin{equation}\label{odeH}
\begin{cases}
\frac{d}{dt}u_i(t) =- \frac 1 {\Delta x} \big( \hat f_{i+\frac12}-\hat f_{i-\frac12}\big),\\
\frac{d}{dt}v_i(t) =- \frac 1 {\Delta x} \big( \hat h_{i+\frac12} -\hat h_{i-\frac12}\big ),
\end{cases}
\end{equation}
where $\hat f_{i+\frac12}$ and $\hat h_{i+\frac12}$ are the numerical fluxes in the governing equation and the derivative equation, respectively.
For the smooth functions $u$ and $v$, the scheme \eqref{odeH} has $k$-th order accuracy if
\begin{equation}\label{semi-kth}
\begin{split}
&f(u)_x|_{x_i} = \frac 1 {\Delta x} \big( \hat f_{i+\frac12}-\hat f_{i-\frac12}\big) + \mathcal{O}(\Delta x)^k,\\
&h(u,v)_x|_{x_i} = \frac 1 {\Delta x} \big( \hat h_{i+\frac12}-\hat h_{i-\frac12}\big) + \mathcal{O}(\Delta x)^{r},\quad r\geq k-1.
\end{split}
\end{equation}
In terms of \cite{js}, the implicit functions $\phi(x)$ and $\psi(x)$ are defined as
\begin{equation*}
\label{impli}
f(u(x))=\frac{1}{\Delta x} \int_{x-\frac{\Delta x}{2}}^{x+\frac{\Delta x}{2}} \phi(\xi) d\xi, \quad  h(u(x),v(x))=\frac 1 {\Delta x} \int_{x-\frac{\Delta x}{2}}^{x+\frac{\Delta x}{2}} \psi(\xi) d\xi.
\end{equation*}
Obviously, we have
\begin{align*}
&f(u)_x|_{x_i} = \frac 1 {\Delta x} \big(\phi(x_{i+\frac12})-\phi(x_{i-\frac12})\big),\quad
h(u,v)_x|_{x_i} = \frac 1 {\Delta x} \big( \psi(x_{i+\frac12})-\psi(x_{i-\frac12})\big).
\end{align*}
This means that the scheme \eqref{odeH} has $k$-th order accuracy if the numerical fluxes
$\hat f_{i+\frac12}$ and $\hat h_{i+\frac12}$ satisfy
\begin{equation}\label{semi-kth-1}
\begin{split}
&\hat f_{i+\frac12} = \phi(x_{i+\frac12}) + \mathcal{O}(\Delta x)^k,\\
&\hat h_{i+\frac12} = \psi(x_{i+\frac12}) + \mathcal{O}(\Delta x)^{r},\quad r\geq k-1.
\end{split}
\end{equation}

We now describe the detailed reconstruction procedure for the numerical fluxes $\hat f_{i+\frac12}$ and $\hat h_{i+\frac12}$ based on $\{u_i,~v_i\}$.
For stability, we should split the fluxes $f(u)$ and $h(u,v)$ into two parts by considering the upwinding mechanism. Given the points value $\{u_i,~v_i\}$, we use the global Lax-Friedrichs splitting:
\begin{equation*}
\begin{split}
&f^{\pm}_i = f^{\pm}(u_i)=\frac12(f(u_i)\pm\alpha u_i) , \quad i=1,...,N_x, \\
&h^{\pm}_i = h^{\pm}(u_i,v_i)=\frac12(h(u_i,v_i)\pm \alpha v_i), \quad i=1,...,N_x,
\end{split}
\end{equation*}
where $\alpha$ is defined as $\max\limits_{u}|f'(u)|$.
The HWENO procedure is applied to $\{f^{\pm}(u), h^{\pm}(u,v)\}$ individually with upwind biased stencils to obtain the numerical fluxes $\hat f^\pm_{i+\frac12}$ and $\hat h^\pm_{i+\frac12}$, and then take $\hat f_{i+\frac12}=\hat f^+_{i+\frac12}+\hat f^-_{i+\frac12}$ and $\hat h_{i+\frac12}=\hat h^+_{i+\frac12}+\hat h^-_{i+\frac12}$.
In this work, the flux $f^{\pm}_{i+\frac12}$ is reconstructed as the convex combination (nonlinear weights) of a quintic polynomial and two quadratic polynomials, while the flux $h^{\pm}_{i+\frac12}$ is reconstructed by the same quintic polynomial directly.

Without loss generality, we here only give the detail reconstruction procedure for $f^{+}_{i+\frac12}$ and $ h^{+}_{i+\frac12}$, while
the procedure for the reconstruction of $f^{-}_{i+\frac12}$ and $ h^{-}_{i+\frac12}$ is the mirror symmetric with respect to $x_{i+\frac12}$.

Following the idea of \cite{ZQ20,ZQd} with artificial linear weights,
we choose a big stencil $T_0 = \{x_{i-1},x_i,x_{i+1}\}$
and two small stencils $T_1 = \{x_{i-1},x_i\}$ and $T_2 = \{x_i,x_{i+1}\}$.
Using the Hermite interpolation on $T_0$ of values $\{u_{i-1},u_{i},u_{i+1},v_{i-1},v_{i},v_{i+1}\}$, there is a quintic polynomial $p_0(x)$ such that
\begin{equation}\label{HWEp0}
\begin{split}
&p_0(x):~\begin{cases}
\frac{1}{\Delta x} \int_{I_{i+\ell}} p_0(x)dx = f^+_{i+\ell},\quad \ell=-1,0,1, \\ \frac{1}{\Delta x} \int_{I_{i+\ell}} p'_0(x)dx = h^+_{i+\ell},\quad  \ell=-1,0,1, \\
\end{cases}
\end{split}
\end{equation}
Similarly, there are two quadratic polynomials $p_1(x)$ and $p_2(x)$ on $T_1$ of values $\{u_{i-1},u_{i},v_{i}\}$ and $T_2$ of values $\{u_{i},u_{i+1},v_{i}\}$, respectively, such that
\begin{equation}\label{HWEp12}
\begin{split}
&p_1(x):~\begin{cases}
\frac{1}{\Delta x} \int_{I_{i+\ell}} p_1(x)dx = f^+_{i+\ell},\quad \ell=-1,0, \\
\frac{1}{\Delta x}\int_{I_{i}} p'_1(x)dx = h^+_{i},
\end{cases}
\\&p_2(x):~\begin{cases}
\frac{1}{\Delta x} \int_{I_{i+\ell}} p_2(x)dx = f^+_{i+\ell},\quad \ell=0,1, \\
\frac{1}{\Delta x}\int_{I_{i}} p'_2(x)dx = h^+_{i}.
\end{cases}
\end{split}
\end{equation}
Evaluate the values of $p_0(x),~p_1(x),~p_2(x)$ and the derivative of $p_0(x)$ at the point $x_{i+\frac12}$, then, we have
\begin{equation*}
\begin{split}
p_0(x_{i+\frac12})&=\frac{11}{60}f_{i-1}+\frac{19}{30}f_i+\frac{11}{60}f_{i+1}+\frac{\Delta x }{20}\big(h_{i-1}+10h_i-h_{i+1}\big),\\
p_1(x_{i+\frac12})&=\frac{1}{6} f_{i-1}+\frac{5}{6}f^+_i
+\frac{2}{3}\Delta xh_i,\\
p_2(x_{i+\frac12})&=\frac{5}{6}f^+_i+\frac{1}{6} f_{i+1}
+\frac{1}{3}\Delta xh_i,
\end{split}
\end{equation*}
and
\begin{equation*}
\begin{split}
p'_0(x_{i+\frac12})&=\frac{1}{4\Delta x}\big(f_{i-1}-8f_i+7f_{i+1}\big)
+\frac{1}{12}\big(h_{i-1}-2h_i-5h_{i+1}\big).
\end{split}
\end{equation*}
The linear weights $\{\gamma_0,~\gamma_1,~\gamma_2\}$ can be chosen as any positive constants with $\gamma_0 +\gamma_1+\gamma_2=1$. Since each of the polynomials in the reconstruction contains $f^+_i$ and $h^+_i$, the HWENO reconstruction could maintain high resolution when the discontinuities occur  at the cell interfaces.
To measure how smooth the functions $p_{\ell}(x),~\ell=0,1,2$ are in the target cell $I_i$,
we compute the smoothness indicators $\beta_{\ell}$ as \cite{js}:
\begin{equation}
\label{GHYZ}
\beta_{\ell}=\sum_{\alpha=1}^k\int_{I_i}{\Delta x}^{2\alpha-1}(\frac{d ^\alpha
p_{\ell}(x)}{d x^\alpha})^2 dx, \quad {\ell}=0,1,2,
\end{equation}
where $k$ is the degree of the polynomials $p_{\ell}(x)$.
The explicit expressions are
\begin{equation}
\label{GHYZP}
\begin{cases}
\beta_0=&\Big(a_1+\frac{1}{4}a_3+\frac{1}{16}a_5\Big)^2
+\frac{13}{3}\Big(a_2+\frac{63}{130}a_4\Big)^2
+\frac{781}{20}\Big(a_3+\frac{8825}{10934}a_5\Big)^2
\\&+\frac{1421461}{2275}a_4^2+\frac{21520059541}{1377684}a_5^2,\\
\beta_1=&(h^+_i\Delta x)^2+\frac{13}{3}\big(\Delta x h^+_{i}-f^+_i+f^+_{i-1}\big)^2,\\
\beta_2=&(h^+_i\Delta x)^2+\frac{13}{3}\big(\Delta x h^+_{i}+f^+_i-f^+_{i-1}\big)^2,\\
\end{cases}
\end{equation}
with
\begin{align*}
\begin{cases}
a_1 = \Delta x
\big(\frac{19}{192} h^+_{i-1}+ \frac{79}{48}h^+_{i}+\frac {19}{192}h^+_{i+1}\big)
+\frac {27}{64}\big(f^+_{i-1}-f^+_{i+1}\big),\\
a_2 = \Delta x
\big(\frac{3}{8}h^+_{i-1}-\frac{3}{8}h^+_{i+1}\big)
+\frac{5}{4}\big(f^+_{i-1}-2f^+_{i}+f^+_{i+1}\big),\\
a_3 = -\Delta x
\big(\frac{11}{24}h^+_{i-1}+\frac{17}{6}h^+_{i}+\frac{11}{24}h^+_{i+1}\big)
+\frac{15}{8}\big(f^+_{i+1}-f^+_{i-1}\big),\\
a_4 = \Delta x\big(\frac{1}{4}h^+_{i+1}-\frac{1}{4}h^+_{i-1}\big)
-\frac{1}{2}\big(f^+_{i-1}-2f^+_{i}+f^+_{i+1}\big),\\
a_5 = \Delta x\big(\frac{1}{4}h^+_{i-1}+h^+_{i}+\frac{1}{4}h^+_{i+1}\big)
+\frac{3}{4}\big(f^+_{i-1}-f^+_{i+1}\big).
\end{cases}
\end{align*}
Following \cite{ZQd}, we define a new parameter $\tau$ to measure the absolute difference
between $\beta_{0}$, $\beta_{1}$ and $\beta_{2}$ as
\begin{equation}\label{tau-fh}
\tau=\frac{1}{4}\Big(|\beta_{0}-\beta_{1}|+|\beta_{0}-\beta_{2}|\Big)^2.
\end{equation}
Then, the nonlinear weights are computed as
\begin{equation*}
\label{90}
\omega_{\ell}=\frac{\bar\omega_{\ell}}{\sum_{\ell=0}^{2}\bar\omega_{\ell}},
\ \mbox{with} \ \bar\omega_{\ell}=\gamma_{\ell}(1+\frac{\tau}{\beta_{\ell}+\varepsilon}), \quad {\ell}=0,1,2.
\end{equation*}
Here, $\varepsilon$ is a small positive number to avoid the denominator by zero. In our computation, we take $\varepsilon = 10^{-6}$ as in the WENO-JS scheme \cite{js} and M-HWENO scheme \cite{ZhaoZhuang-2020JSC-FD}. Finally, the values of $\hat{f}^+_{i+\frac12}$ and $\hat h^+_{i+\frac12}$ are reconstructed by
\begin{align*}
\begin{cases}
\hat{f}^+_{i+\frac12} =\omega_0 \Big( \frac 1 {\gamma_0}p_0(x_{i+\frac12})  -
\frac {\gamma_{1}} {\gamma_0} p_{1}(x_{i+\frac12})
-\frac {\gamma_{2}} {\gamma_0} p_{2}(x_{i+\frac12}) \Big)
+\omega_{1} p_{1}(x_{i+\frac12})+\omega_{2} p_{2}(x_{i+\frac12}),\\
\hat{h}^+_{i+\frac12} =p'_0(x_{i+\frac12}).
\end{cases}
\end{align*}
Obviously, we have $|\hat{f}_{i+\frac12} - \phi(x_{i+\frac12})| = \mathcal{O}(\Delta x)^6$ and $|\hat{h}^+_{i+\frac12} - \psi(x_{i+\frac12})| = \mathcal{O}(\Delta x)^5$.
From \eqref{semi-kth} and \eqref{semi-kth-1}, the semi-discrete scheme \eqref{odeH} at least has the fifth-order accuracy for the smooth functions $u$ and $v$.

\subsection{HWENO limiter for the solution derivative}
\label{sec:HWENO-limiter}

Since the solution for hyperbolic conservation laws often contains discontinuities, the derivative of the solution would be quite large near discontinuities, then, it is a natural idea that we should deal with the derivative values carefully. Several works have been done to control it in the finite difference HWENO framework.
For example, Liu and Qiu  \cite{LiuQ1} (the first  finite difference HWENO scheme) used the different polynomials in the reconstruction to escape discontinuities.
Zhao et al. \cite{ZhaoZhuang-2020JSC-FD} (M-HWENO scheme) modified the derivatives before the reconstruction.
Li et al. used the center point value to reconstruct the fluxes automatically near discontinuities in the multi-resolution HWENO scheme \cite{LiMRHW1}.
However, the schemes \cite{LiMRHW1,LiuQ1} only achieve the fourth-order accuracy in two dimensions.

To both avoid spurious oscillations and maintain the fifth-order accuracy,  we add an HWENO limiter to control the derivatives following the idea of \cite{ZhaoZhuang-2020JSC-FD}.
Instead of using the modified derivatives both in fluxes reconstruction and time discretization as in \cite{ZhaoZhuang-2020JSC-FD}, we only apply the modified derivatives in time discretization while remaining the original derivatives in fluxes reconstruction.
It is interesting that $h^+_{i+\frac12}$ can be approximated by a quintic polynomial directly in the proposed HWENO scheme, while $h^+_{i+\frac12}$ must be reconstructed by a nonlinear HWENO method in M-HWENO scheme \cite{ZhaoZhuang-2020JSC-FD}.
The HWENO limiter for the derivative is based on the convex combination of a quartic polynomial with two linear polynomials, and the linear weights also can be chosen as artificial positive number as long as their sum equals one.
Now, we describe the detail of the HWENO limiter to control $v_i$ and obtain the modified derivative $\tilde{v}_i$ finally.
Using the Hermite interpolation on stencils $T_0$, $T_1$ and $T_2$, respectively, there are a unique quartic polynomial $q_0(x)$ and two linear polynomials $q_1(x)$ and $q_2(x)$, such that
\begin{equation*}
\begin{split}
&q_0(x): ~\begin{cases}
q_0(x_{i+\ell})=u_{i+\ell},  \quad \ell=-1,0,1,\\
q'_0(x_{i+\ell})= v_{i+\ell},\quad\ell=-1,1,\\
\end{cases}\\
&q_1(x): ~q_1(x_{i+\ell})=u_{i+\ell}, \quad \ell=-1,0, \\
&q_2(x): ~q_2(x_{i+\ell})=u_{i+\ell}, \quad \ell=0,1.
\end{split}
\end{equation*}
And then, we have
\begin{equation*}
\begin{split}
q'_0(x_i)&=
\frac{3}{4 \Delta x}\big(u_{i+1}-u_{i-1}\big)-\frac{1}{4}\big(v_{i-1}+v_{i+1}\big),\\
q'_1(x_i)&=\frac{1}{\Delta x}\big(u_i-u_{i-1}\big),\\
q'_2(x_i)&=\frac{1}{\Delta x}\big(u_{i+1}-u_{i}\big),\\
\end{split}
\end{equation*}
where the linear weights $d_0$, $d_1$, $d_2$ can be chosen as any positive constants with $d_0 +d_1+d_2=1$.

Similarly as described in Section \ref{sec:HWENO-1d-u}, we compute the smoothness indicators $\beta_{\ell}$ to
measure how smooth the functions $q_{\ell}(x),~\ell=0,1,2$
are in the target cell $I_i$ as:
\begin{equation}\label{GHYZ-q}
\beta_{\ell}=\sum_{\alpha=1}^{r}\int_{I_i}{\Delta x}^{2\alpha-1}(\frac{d ^\alpha
q_{\ell}(x)}{d x^\alpha})^2 dx, \quad {\ell}=0,1,2,
\end{equation}
where $r$ is the degree of the polynomials $q_{\ell}(x)$.
The explicit formulas are given by
\begin{align*}
\beta_0&=\big(a_1+\frac{1}{4}a_3\big)^2
+\frac{13}{3}\big(a_2+\frac{63}{130}a_4\big)^2
+\frac{781}{20}a_3^2+\frac{1421461}{2275}a_4^2,\\
\beta_1&=(u_{i}-u_{i-1})^2, \\
\beta_2&=(u_i -u_{i+1})^2,
\end{align*}
with
\begin{align*}
\begin{cases}
a_1 =
-\frac{\Delta x}{4} \big( v_{i-1}+v_{i+1}\big)
+\frac{3}{4}\big(u_{i+1}-u_{i-1}\big),\\
a_2 =
\frac{\Delta x}{4} \big(v_{i-1}-v_{i+1}\big)
+u_{i-1}-2u_{i}+u_{i+1},\\
a_3 =
\frac{\Delta x}{4} \big(v_{i-1}+v_{i+1}\big)
+\frac{1}{4}\big(u_{i-1}-u_{i+1}\big),\\
a_4 = \frac{\Delta x}{4} \big(v_{i+1}-v_{i-1}\big)
-\frac{1}{2}\big(u_{i-1}-2u_{i}+u_{i+1}\big).
\end{cases}
\end{align*}
The nonlinear weights are defined as
\begin{align*}\label{9}
&\lambda_\ell=\frac{\bar\lambda_\ell}{\sum_{\ell=0}^{2}\bar\lambda_{\ell}},\qquad
\bar\lambda_{\ell}=d_{\ell}\Big(1+\frac{\tau}{\beta_{\ell}+\varepsilon}\Big),\quad \ell=0,1,2,
\quad
\end{align*}
where $\tau = \frac{1}{4}\big(|\beta_{0}-\beta_{1}|+|\beta_{0}-\beta_{2}|\big)^2$
and $\varepsilon = 10^{-6}$ is to avoid the denominator by zero.
Finally, the modified derivative $ \tilde{v}_i$ is defined as
\begin{equation*}
 \tilde{v}_i = \lambda_0 \Big(
 \frac{1}{d_0}q'_0(x_i) - \frac{d_1}{d_0} q'_{1}(x_i)-\frac{d_2}{d_0} q'_{2}(x_i) \Big)
 + \lambda_1 q'_1(x_i)+ \lambda_2 q'_2(x_i).
\end{equation*}
Obviously, we have $|\tilde{v}_i - v_i| = \mathcal{O}(\Delta x)^4$. From \eqref{semi-kth} and \eqref{semi-kth-1}, it is not difficult to know that it maintains the fifth-order accuracy of the HWENO scheme.

\vspace{10pt}

Denote
\begin{equation}
\begin{split}
\mathcal{L}^1_{i}(u,v) = - \frac 1 {\Delta x} \big( \hat f_{i+\frac12}-\hat f_{i-\frac12}\big),\quad i=1,...,N_x,\\
\mathcal{L}^2_{i}(u,v) = - \frac 1 {\Delta x} \big( \hat h_{i+\frac12}-\hat h_{i-\frac12}\big),\quad i=1,...,N_x.
\end{split}
\end{equation}
For time discretization of \eqref{odeH}, we use the explicit third-order SSP Runge-Kutta scheme, then we have the fully-discrete scheme as, for $i=1,...,N_x$
\begin{subequations}
\begin{align}
&\begin{cases}
u^{(1)}_i =u^n_i+ \Delta t \mathcal{L}^1_{i}(u^{n},v^{n}),\\
v^{(1)}_i =\tilde{v}_i ~+ \Delta t \mathcal{L}^2_{i}(u^{n},v^{n}),
\end{cases}\\
&\begin{cases}
u^{(2)}_i =\frac{3}{4}u^n_i+
\frac{1}{4}\big( u^{(1)}_i+ \Delta t \mathcal{L}^1_{i}(u^{(1)},v^{(1)})\big),\\
v^{(2)}_i =\frac{3}{4}\tilde{v}_i ~+
\frac{1}{4}\big( \tilde{v}^{(1)}_i ~+ \Delta t  \mathcal{L}^2_{i}(u^{(1)},v^{(1)})\big),
\end{cases}\\
&\begin{cases}
u^{n+1}_i =\frac{1}{3}u^n_i+
\frac{2}{3}\big( u^{(2)}_i+ \Delta t \mathcal{L}^1_{i}(u^{(2)},v^{(2)})\big),\\
v^{n+1}_i =\frac{1}{3}\tilde{v}_i ~+
\frac{2}{3}\big( \tilde{v}^{(2)}_i~ + \Delta t  \mathcal{L}^2_{i}(u^{(2)},v^{(2)})\big),
\end{cases}
\end{align}
\end{subequations}
where $\tilde{v}_i$ is the modified derivative of $v_i$ obtained by the HWENO limiter.

\begin{rem}
For the system case, such as the one-dimensional compressible Euler equations, the HWENO reconstruction for $f^+_{i+\frac12}$ is implemented  based on the local characteristic decomposition \cite{js},  while the linear approximation for $h^+_{i+\frac12}$ is performed by  component-by-component.
\end{rem}


\vspace{10pt}

A major advantage of the high-order finite difference scheme is that it is straightforward to extend the scheme in one-dimension to two-dimensions by dimension-by-dimension.
Hence, we also can extend the proposed finite difference HWENO scheme to two-dimensions easily, but use the special treatments of the mixed derivative terms as in \cite{ZhaoZhuang-2020JSC-FD}.
 One striking feature of the proposed HWENO scheme is that it can achieve the fifth-order accuracy in two dimensions, while other finite difference HWENO schemes, e.g., \cite{LiMRHW1, LiuQ1}, only have the fourth-order accuracy.

We consider the two-dimensional scalar hyperbolic conservation laws
\begin{equation}\label{EQ2}
\begin{cases}
u_t+ f(u)_x+g(u)_y=0, \\
u(x,y,0)=u_0(x,y).
\end{cases}
\end{equation}
We rewrite \eqref{EQ2} by bringing its derivative equations as
\begin{equation}\label{EQ3}
\begin{cases}
u_t+ f(u)_x+g(u)_y=0, \\
v_t+h(u,v)_x+\xi(u,v)_y=0, \\
w_t+\eta(u,w)_x+\theta(u,w)_y=0, \
\end{cases}
\end{equation}
where
\begin{align*}
&v=u_x,\quad h(u,v)=f'(u)u_x =f'(u)v, \quad \xi(u,v)=g'(u)u_x = g'(u)v,
\\&w=u_y,\quad \eta(u,w)=f'(u)u_y =f'(u)w, \quad \theta(u,w)=g'(u)u_y =g'(u)w.
\end{align*}
Denote $I_{i} = [x_{i-\frac12},x_{i+\frac12}]$, $J_{j}=[y_{j-\frac12},y_{j+\frac12}]$, $I_{i,j} = I_{i}\times J_j $, $\Delta x=x_{i+\frac12}-x_{i-\frac12}$, $\Delta y=y_{j+\frac12}-y_{j-\frac12}$, and $(x_i,y_j)$ is the center of the element $I_{i,j}$.
The semi-discrete finite difference scheme of \eqref{EQ3} is
\begin{equation}
\label{ode2H}
  \begin{cases}
   \frac{d}{dt}u_{i,j}(t) &=- \frac 1 {\Delta x} \Big( \hat f_{i+\frac12,j}-\hat f_{i-\frac12,j}\Big)
   -\frac1{\Delta y}\Big( \hat g_{i,j+\frac12}-\hat g_{i,j-\frac12}\Big),\\
   \frac{d}{dt}v_{i,j}(t) &=- \frac 1 {\Delta x} \Big( \hat h_{i+\frac12,j}-\hat h_{i-\frac12,j}\Big)
   -\frac1{\Delta y}\Big( \hat \xi_{i,j+\frac12}-\hat \xi_{i,j-\frac12}\Big),\\
   \frac{d}{dt}w_{i,j}(t) &=- \frac 1 {\Delta x} \Big( \hat \eta_{i+\frac12,j}-\hat \eta_{i-\frac12,j}\Big)
   -\frac1{\Delta y}\Big( \hat \theta_{i,j+\frac12}-\hat \theta_{i,j-\frac12}\Big).
  \end{cases}
\end{equation}
Here, the numerical fluxes $ \hat f_{i\pm\frac12,j}$, $\hat g_{i,j\pm\frac12}$, $\hat h_{i\pm\frac12,j}$ and $\hat \theta_{i,j\pm\frac12}$ are reconstructed by a dimension-by-dimension manner, seen in Section \ref{sec:HWENO-1d-u}. We can get
\begin{equation}
\begin{split}\label{dd-app}
&f(u)_x|_{(x_{i},y_j)} = \frac 1 {\Delta x} \Big( \hat f_{i+\frac12,j}-\hat f_{i-\frac12,j}\Big) + \mathcal{O}(\Delta x^6),\\
&g(u)_y|_{(x_{i},y_j)} = \frac 1 {\Delta y} \Big( \hat g_{i,j+\frac12}-\hat g_{i,j-\frac12}\Big) + \mathcal{O}(\Delta y^6),\\
&h(u,v)_x|_{(x_{i},y_j)}=\big(f'(u)u_x\big)_x|_{(x_{i},y_j)} = \frac 1 {\Delta x} \Big( \hat h_{i+\frac12,j}-\hat h_{i-\frac12,j}\Big)+ \mathcal{O}(\Delta x^5),\\
&\theta(u,w)_y|_{(x_{i},y_j)}=\big(g'(u)u_y\big)_y|_{(x_{i},y_j)} = \frac1{\Delta y}\Big( \hat \theta_{i,j+\frac12}-\hat \theta_{i,j-\frac12}\Big)+ \mathcal{O}(\Delta y^5).
\end{split}
\end{equation}
To ensure the fifth-order accuracy of the scheme \eqref{ode2H}, we would like to find at least fourth-order approximations for the mixed derivative terms $\xi(u,v)_y$ and $\eta(u,w)_x$ at point $(x_{i},y_j)$.
However, $\xi(u,v)_y$ and $\eta(u,w)_x$ don't have their primitive functions in $y$ and $x$ directions, respectively, thus, they can't be approximated using the same way as was down in other derivative terms
$h_(u,v)_x$ and $\theta_(u,w)_y$. Here, we adopt the linear approximation directly (without fluxes splitting) for the mixed derivative terms as
\begin{equation}\label{mixed}
\begin{split}
&\hat \xi_{i,j+\frac12}=
-\frac{1}{12}\xi_{i,j-1}+\frac{7}{12}\xi_{i,j}
+\frac{7}{12}\xi_{i,j+1}-\frac{1}{12}\xi_{i,j+2},
\\&\hat \eta_{i+\frac12,j}=
-\frac{1}{12}\eta_{i-1,j}+\frac{7}{12}\eta_{i,j}
+\frac{7}{12}\eta_{i+1,j}-\frac{1}{12}\eta_{i+2,j}.
\end{split}
\end{equation}
Then, we have
\begin{equation}
\begin{split}\label{md-app}
&\xi(u,v)_y|_{(x_{i},y_j)}=\big(g'(u)u_x\big)_y|_{(x_{i},y_j)} = \frac 1 {\Delta x} \Big( \hat \xi_{i,j+\frac12}-\hat \xi_{i,j-\frac12}\Big)+ \mathcal{O}(\Delta y^4),\\
&\eta(u,w)_x|_{(x_{i},y_j)}=\big(f'(u)u_y\big)_x|_{(x_{i},y_j)} = \frac1{\Delta y}\Big( \hat \eta_{i+\frac12,j}-\hat \eta_{i-\frac12,j}\Big)+ \mathcal{O}(\Delta x^4).
\end{split}
\end{equation}
Thus, from \eqref{semi-kth} and \eqref{semi-kth-1}, we can prove that the HWENO scheme \eqref{ode2H} has fifth-order accuracy.

Similarly as in one dimension, we also use the explicit third-order SSP Runge-Kutta scheme to discretize \eqref{odeH}, and add the HWENO limiter to control the derivatives $v$ and $w$ in time discretization by a dimension-by-dimension manner (cf. Section \ref{sec:HWENO-limiter}).

\begin{rem}
For the system case, the HWENO  procedures in $x$ and $y$ directions are implemented on each local characteristic direction, respectively, while the linear approximations for the fluxes in the derivative equations are performed on each component straightforwardly.
\end{rem}

\section{Numerical experiments}
\label{sec:numerical}

In this section, we present the numerical results to show the good performances of the proposed finite difference HWENO scheme combined with limiter.
For comparisons, we consider three variants of the HWENO or WENO schemes:
\begin{itemize}
  \item The fifth-order finite difference L-HWENO scheme: the proposed HWENO scheme, where the modified derivatives are only used in time discretization while remaining the original derivatives in fluxes reconstruction, seen in Section \ref{sec:HWENO}.
  \item The fifth-order finite difference M-HWENO scheme: the modified HWENO scheme \cite{ZhaoZhuang-2020JSC-FD}, where the modified derivatives are used both in fluxes reconstruction and time discretization.
  \item The fifth-order finite difference WENO-JS scheme: the classical WENO scheme constructed  by Jiang and Shu \cite{js}.
\end{itemize}
Since the M-HWENO scheme \cite{ZhaoZhuang-2020JSC-FD} had shown its smaller error for accuracy tests and higher resolution near discontinuities than the WENO-JS scheme \cite{js}, we only present the results of the WENO-JS scheme for the accuracy tests to compare their efficiency here.

The CFL number is set as $0.6$ in our computation. Unless otherwise stated, the linear weights for the reconstruction fluxes in the governing equation and limiter are taken as $\{\gamma_0=0.98,~\gamma_1=\gamma_2=0.01 \}$ and $\{d_0=0.98,~d_1=d_2=0.01\}$, respectively.
For examples where the analytical exact solution is unavailable, we take the numerical solution obtained by the WENO-JS scheme \cite{js} with $N_x=2000$ as the  referenced ``exact'' solution.


\begin{example} \label{burgers-1d}
(Accuracy test of the one-dimensional Burgers' equation.)
\end{example}
This example is used to verify the fifth-order accuracy and efficiency of the proposed L-HWENO scheme for the one-dimensional nonlinear Burgers' equation.
The Burgers' equation in one dimension reads as
\begin{equation}\label{1dbugers}
  u_t+(\frac {u^2} 2)_x=0, \quad 0<x<2,
\end{equation}
with the periodic boundary condition. The initial condition is $u(x,0)=0.5+\sin(\pi x)$.

The final simulation time is $T=0.5/\pi$ when the solution is still smooth.
The $L^1$ and $L^\infty$ norm of the error with M-HWENO and L-HWENO schemes are listed in Table \ref{tburgers1d}.
It can be seen that the schemes both achieve the optimal fifth-order accuracy, but  the error obtained by the L-HWENO scheme is smaller than that obtained by M-HWENO scheme.

\begin{figure}[H]
\centering
\subfigure[Example \ref{burgers-1d}.]{
\includegraphics[width=0.45\textwidth,trim=20 10 30 30,clip]{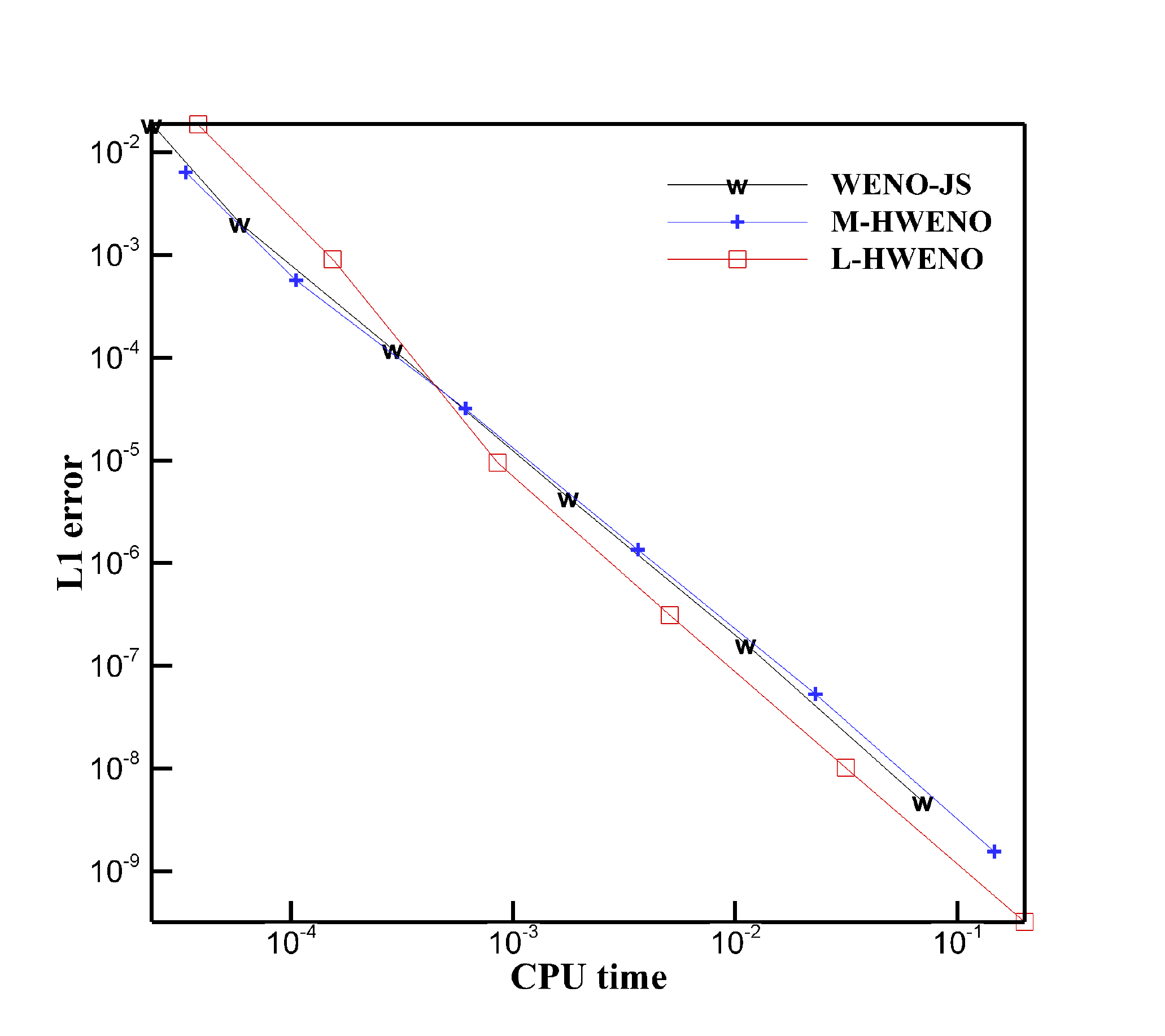}}
\subfigure[Example \ref{Euler-1d}.]{
\includegraphics[width=0.45\textwidth,trim=20 10 30 30,clip]{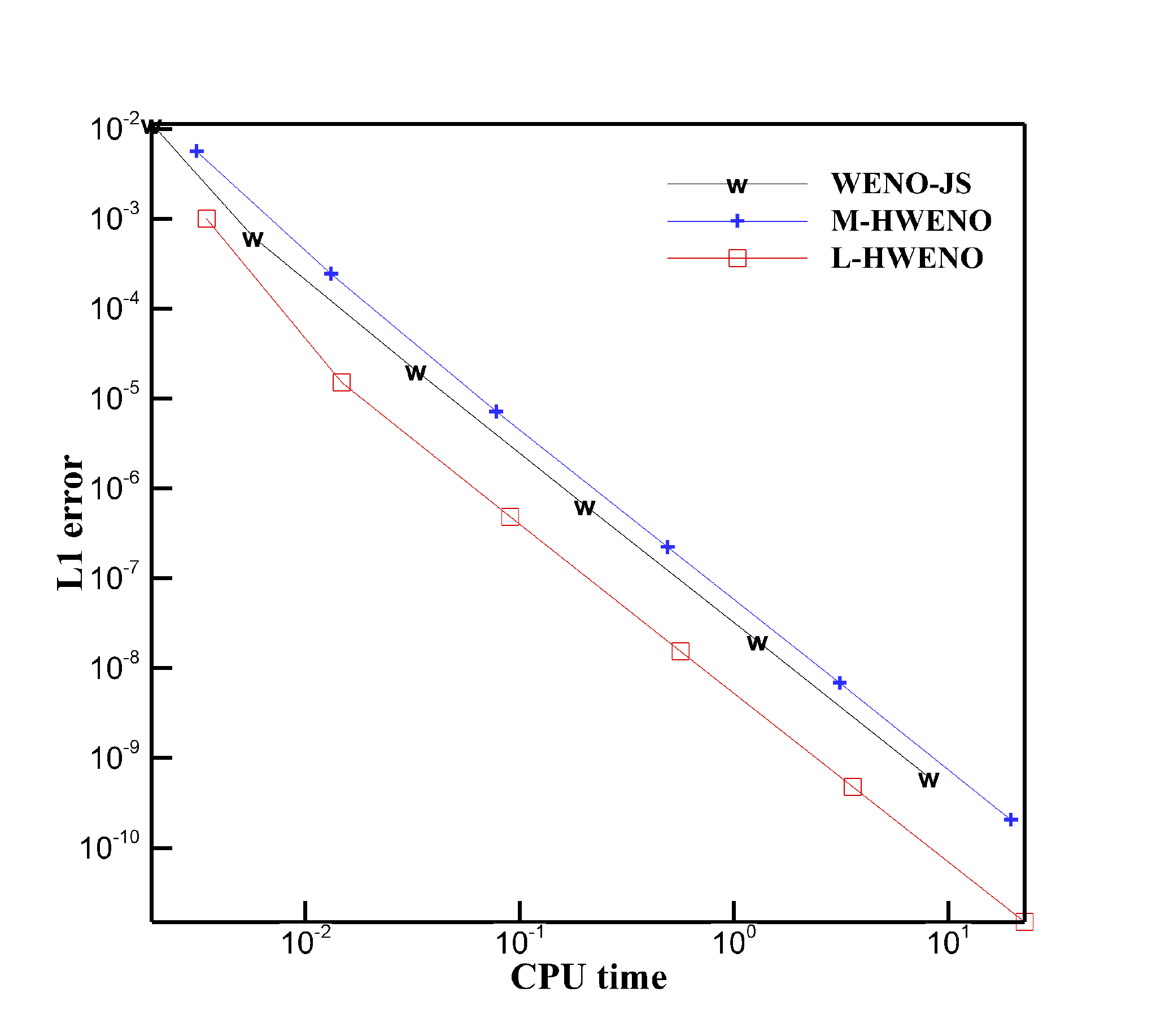}}
 \caption{The error of $L^1$ norm against the CPU time.}
\label{1d-CPU}
\end{figure}

For comparisons of their efficiency, we also plot the $L^1$ norm of the error against CPU time measured in seconds in Fig.~\ref{1d-CPU}(a).
One can find that the proposed L-HWENO scheme is more efficient than the M-HWENO scheme or WENO-JS  scheme in the sense that the former leads to a smaller error than the latter for a fixed amount of the CPU time. Thus, the proposed L-HWENO scheme has better efficiency than the M-HWENO scheme when $N_x$ is relatively large.


\begin{table}[H]
\caption{Example \ref{burgers-1d}. The $L^1$ and $L^\infty$ norm of the error computed by the M-HWENO and L-HWENO schemes.}
\centering
\vspace{5pt}
\begin{tabular} {ccccccccc}
\toprule
 $N_x$& \multicolumn{4}{c}{M-HWENO}&\multicolumn{4}{c}{L-HWENO}\\
\cline{2-5} \cline{6-9}
 & $L^1$ error &  order & $L^\infty$ error & order &
  $L^1$ error &  order & $L^\infty$ error  &order  \\
\midrule
     10 &     6.42E-03 &            &     3.01E-02 &
        &     1.88E-02 &            &     6.17E-02 &      \\
     20	&     5.69E-04 &       3.50 &     4.05E-03 &       2.90
        &     9.18E-04 &       4.36 &     5.42E-03 &       3.51 \\
     40 &     3.20E-05 &       4.15 &     3.86E-04 &       3.39
        &     9.52E-06 &       6.59 &     1.13E-04 &       5.59 \\
     80 &     1.35E-06 &       4.56 &     1.31E-05 &       4.88
        &     3.13E-07 &       4.93 &     4.14E-06 &       4.77 \\
    160 &     5.32E-08 &       4.67 &     5.67E-07 &       4.53
        &     1.02E-08 &       4.94 &     1.34E-07 &       4.95 \\
    320 &     1.55E-09 &       5.10 &     2.06E-08 &       4.78
        &     3.23E-10 &       4.98 &     4.29E-09 &       4.96 \\
\bottomrule
\end{tabular}
\label{tburgers1d}
\end{table}

\begin{example} \label{Euler-1d}
(Accuracy test of the one-dimensional Euler equations.)
\end{example}
This example is used to verify the fifth-order accuracy and efficiency of the proposed L-HWENO scheme for the one-dimensional system of Euler equations.
The Euler equations in one-dimension read as
\begin{equation}\label{euler1}
 \frac{\partial}{\partial t}
 \begin{bmatrix*}[c]
 \rho \\
 \rho \mu \\
 E
 \end{bmatrix*}
 +
 \frac{\partial}{\partial x}
\begin{bmatrix*}[c]
\rho \mu \\
\rho \mu^{2}+p \\
\mu(E+p)
\end{bmatrix*}
=0,
\end{equation}
where $\rho$ is the density, $\mu$ is the velocity, $E$ is the total energy and $p$ is the pressure, where $E = p/(\gamma-1)+\rho\mu^2/2 $ with $\gamma=1.4$. The computational domain is $[0, 2]$. The periodic boundary conditions are used for all unknown variables and the initial conditions are given by
\[
\rho(x,0)=1+0.2\sin(\pi x),\quad \mu(x,0)=1,\quad p(x,0)=1.\]
The exact solution of this example is
\[\rho(x,t)=1+0.2\sin(\pi(x-t)),\quad \mu(x,0)=1,\quad p(x,0)=1.\]
The final simulation time is $T=2$.

The $L^1$ and $L^\infty$ norm of the error obtained by the M-HWENO and the proposed L-HWENO schemes are presented in Table \ref{tEluer1d}.
Similarly as the last example, we can clearly see that the schemes both achieve the optimal fifth-order accuracy, and the error of the solution obtained by the L-HWENO scheme is smaller than that obtained by the M-HWENO scheme.

To show the efficiency of the proposed L-HWENO scheme for the one-dimensional system,
we plot the $L^1$ norm of the error against CPU time in Fig.~\ref{1d-CPU}(b).
One can find that the L-HWENO scheme is more efficient than either the M-HWENO or  WENO-JS  scheme in the sense that the former leads to a smaller error than the latter for a fixed amount of the CPU time. Moreover, the better efficiency of the proposed L-HWENO scheme is more obvious in this example than that in Example \ref{burgers-1d} for one-dimensional Burgers' equation. The reason is that the linear approximations for the fluxes in the L-HWENO scheme are applied in each component directly without any local characteristic decomposition for systems.

\begin{table}[H]
\caption{Example \ref{Euler-1d}. The $L^1$ and $L^\infty$ norm of the error obtained by the M-HWENO and L-HWENO schemes.}
\centering
\vspace{5pt}
\begin{tabular} {ccccccccc}
\toprule
$N_x$ & \multicolumn{4}{c}{M-HWENO}& \multicolumn{4}{c}{L-HWENO}\\
\cline{2-5} \cline{6-9}
  &$L^1$ error &  order & $L^\infty$ error &  order &
  $L^1$ error &  order & $L^\infty$ error  & order \\
  \midrule
10 &  5.62E-03 &       & 9.93E-03 &
   &  9.97E-04 &       & 2.49E-03 &\\
20 &  2.43E-04 & 4.53  & 4.42E-04 & 4.49
   &  1.51E-05 & 6.05 & 3.51E-05 & 6.15 \\
40 & 7.14E-06 & 5.09  & 1.40E-05 & 4.98
   & 4.83E-07 & 4.96  & 8.35E-07 & 5.39\\
80 & 2.21E-07 & 5.01  & 4.36E-07 & 5.00
   & 1.53E-08 & 4.98  & 2.46E-08 & 5.08\\
160& 6.84E-09 & 5.01  & 1.25E-08 & 5.13
   & 4.80E-10 & 4.99  & 7.59E-10 & 5.02 \\
320& 2.05E-10 & 5.06  & 3.53E-10 & 5.14
   & 1.51E-11 & 4.99  & 2.37E-11 & 5.00\\
   \bottomrule
\end{tabular}
\label{tEluer1d}
\end{table}

\begin{figure}[H]
\subfigure[solution $u$: $N_x=80$]{
\includegraphics[width=0.45\textwidth,trim=20 10 30 30,clip]{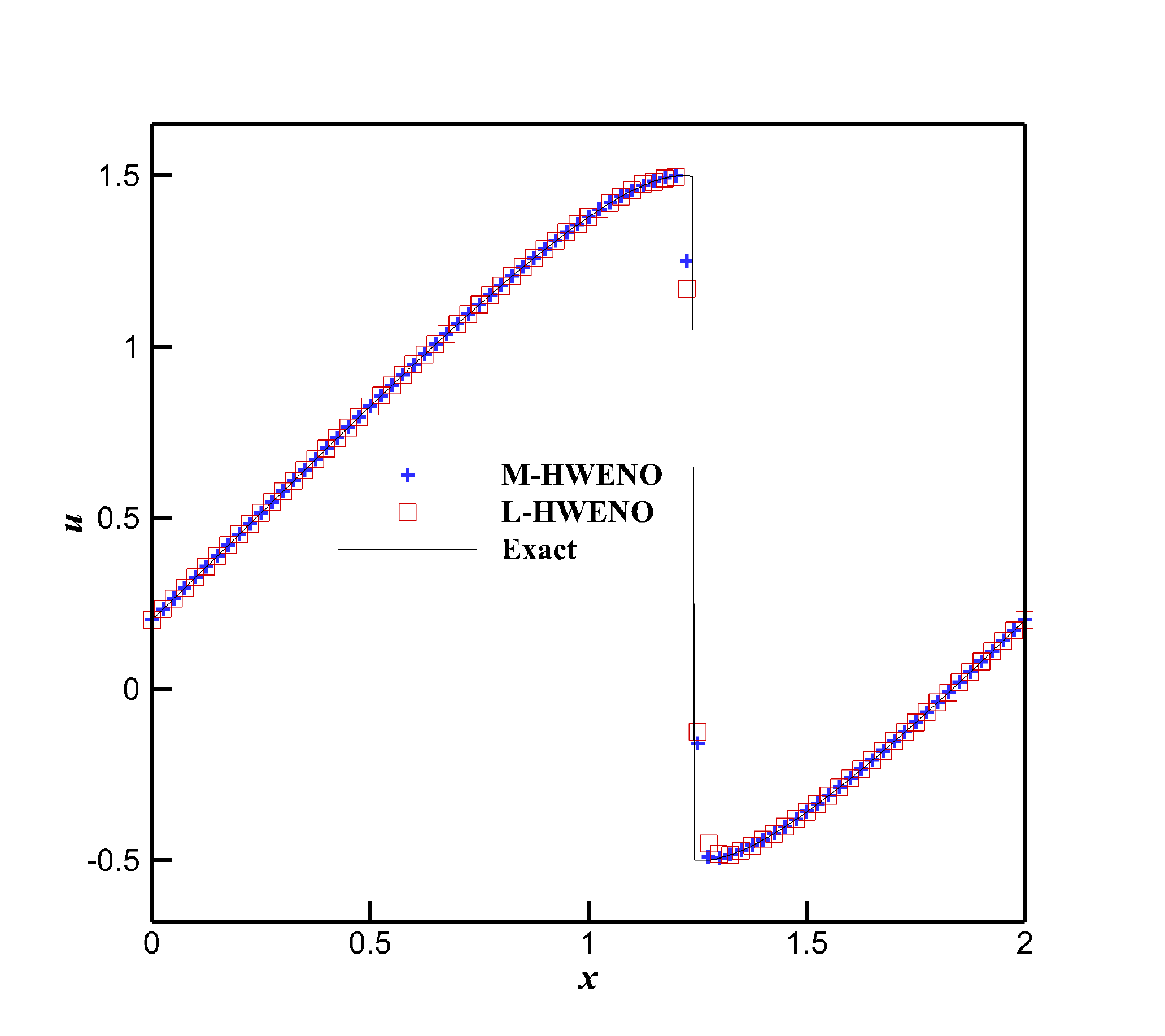}}
\subfigure[close view of (a)]{
\includegraphics[width=0.45\textwidth,trim=20 10 30 30,clip]{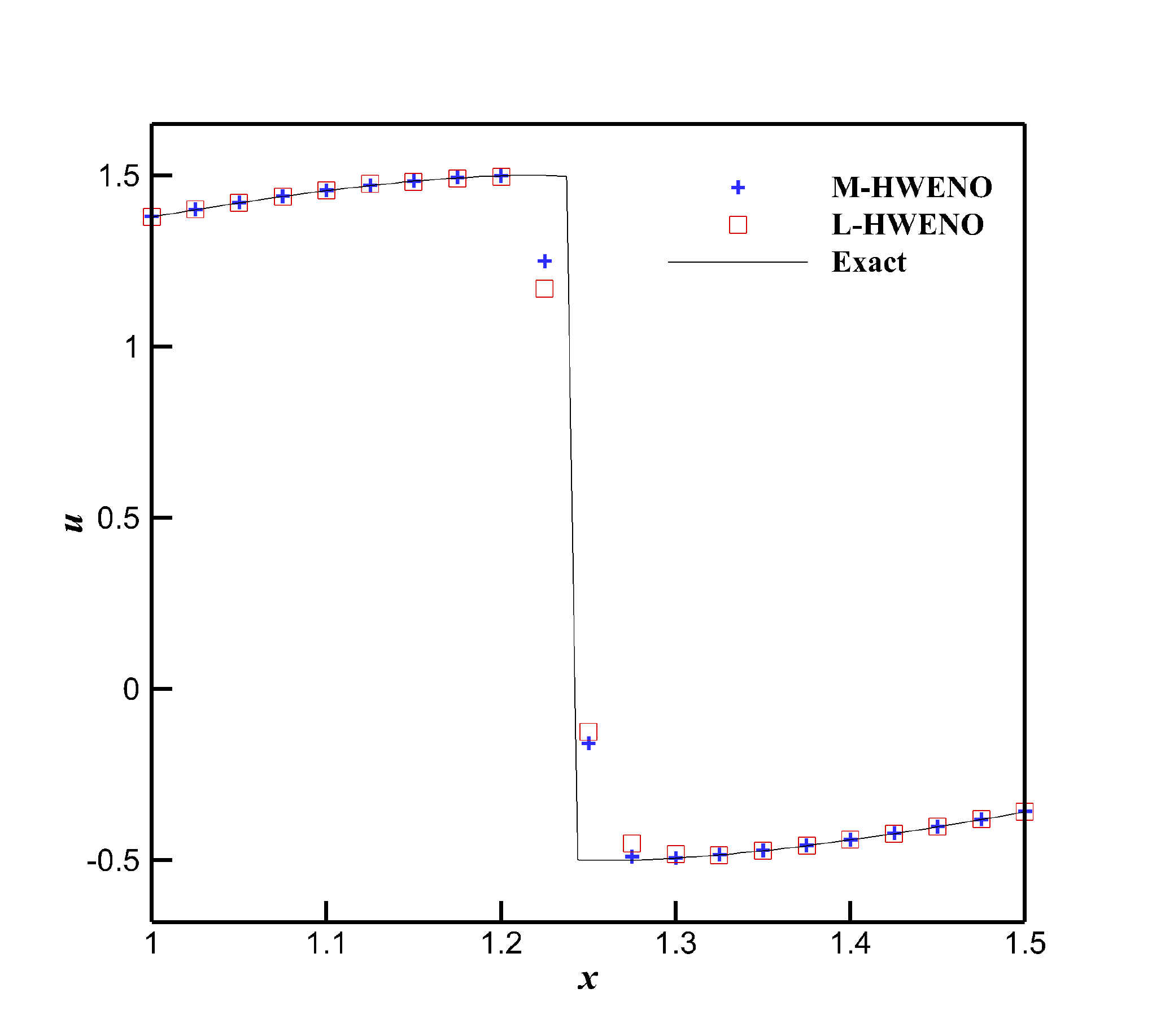}}
 \caption{Example~\ref{burgers-1d-shock}. The numerical solution at $T=1.5/\pi$ obtained by the M-HWENO and L-HWENO schemes.}
\label{Fburges1d}
\end{figure}

%
\begin{example} \label{burgers-1d-shock}
(Shock wave of the one-dimensional Burgers' equation.)
\end{example}
In this test, we repeat the one-dimensional Burgers' equation \eqref{1dbugers} given in Example \ref{burgers-1d}, but the final simulation time is $T=1.5/\pi$ when the solution is discontinuous.
The numerical solution $u$ obtained by the M-HWENO and L-HWENO schemes against the exact solution is plotted in Fig.~\ref{Fburges1d}. From the figure, we can know that their performances are similar with non-oscillations.

\begin{example} \label{BL-1d}
(Buckley-Leverett problem of the one-dimensional nonlinear non-convex equation.)
\end{example}
This example is used to verify the performance of the proposed L-HWENO scheme for the one-dimensional nonlinear non-convex scalar equation. It is not easy to simulate since the numerical solution may violate the entropy condition.
We consider a one-dimensional nonlinear non-convex Buckley-Leverett problem
\begin{equation*}
  u_t+\left( \frac{4u^2}{4u^2+(1-u)^2}\right)_x=0,\quad -1\leq x\leq 1.
\end{equation*}
The initial condition is
\begin{equation*}
u(x,0) =
\begin{cases}
1, \quad \hbox{for~} -\frac12\leq x\leq0,\\
0, \quad \hbox{otherwise.}
\end{cases}
\end{equation*}
The final simulation time is $T=0.4$. The exact solution of this problem contains both shock wave and rarefaction wave.

\begin{figure}[H]
\subfigure[solution $u$: $N_x=80$]{
\includegraphics[width=0.45\textwidth,trim=20 10 30 30,clip]{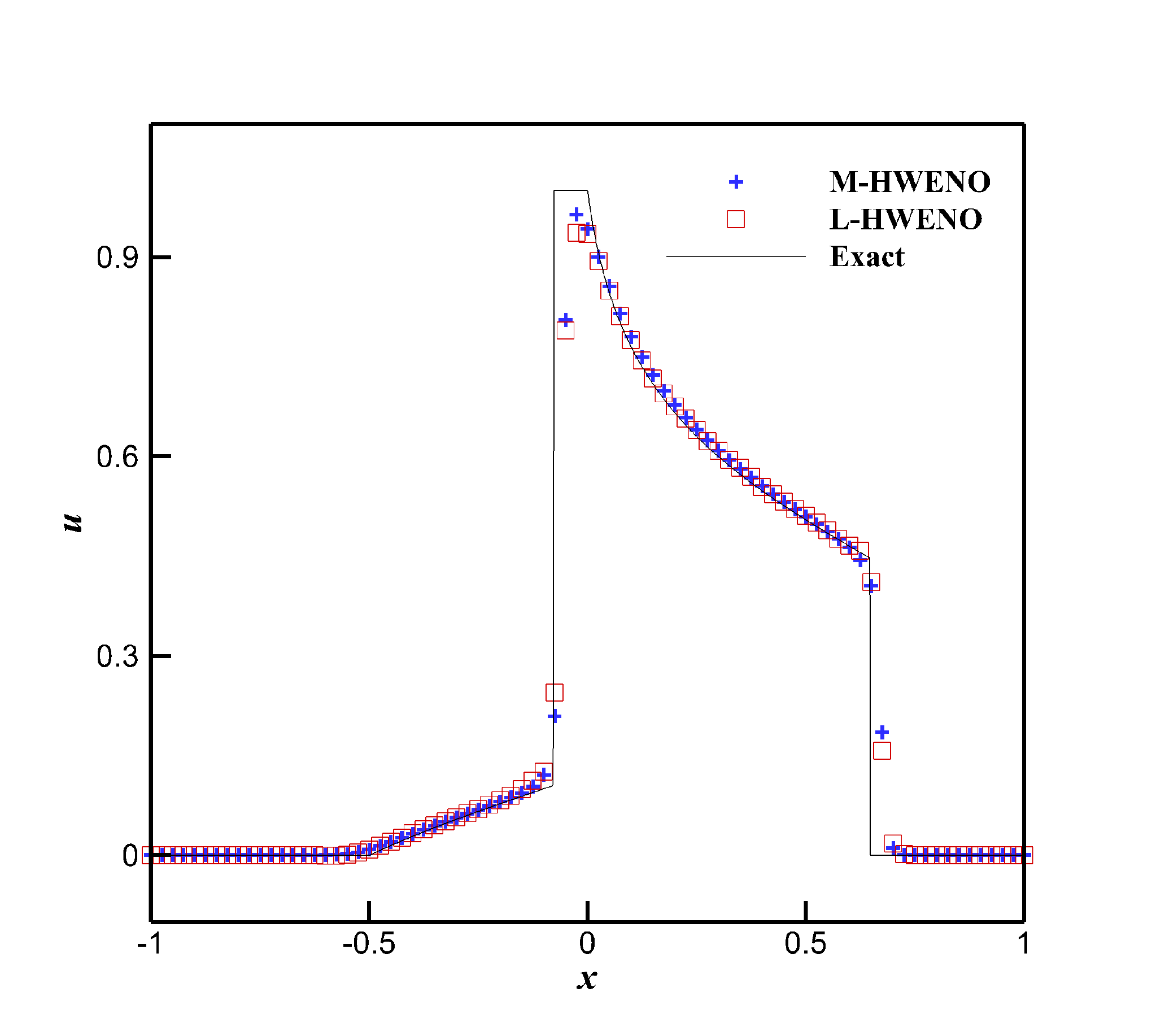}}
\subfigure[close view of (a)]{
\includegraphics[width=0.45\textwidth,trim=20 10 30 30,clip]{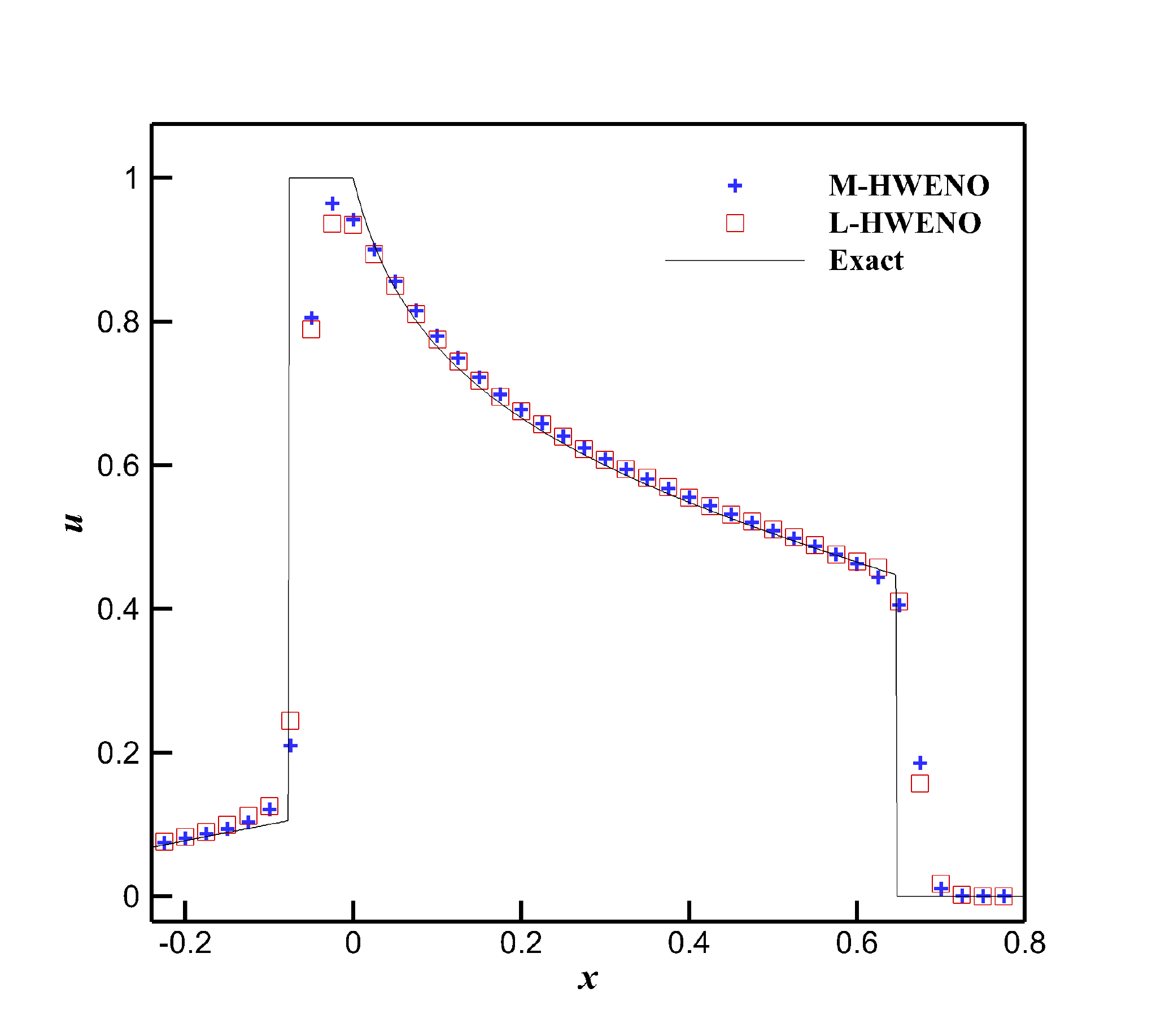}}
 \caption{Example~\ref{BL-1d}. The numerical solution at $T=0.4$ obtained by the M-HWENO and L-HWENO schemes.}
\label{FBuk_levert}
\end{figure}

In Fig.~\ref{FBuk_levert}, we plot the solution obtained by the M-HWENO and L-HWENO schemes. We can see that the two schemes have similar performances generally, but the M-HWENO scheme seemly has higher resolution near the peak, and we will investigate the reasons of this phenomenon in Example \ref{blast-1d} below.

\begin{figure}[H]
\subfigure[density $\rho$: $N_x=200$]{
\includegraphics[width=0.45\textwidth,trim=20 10 30 30,clip]{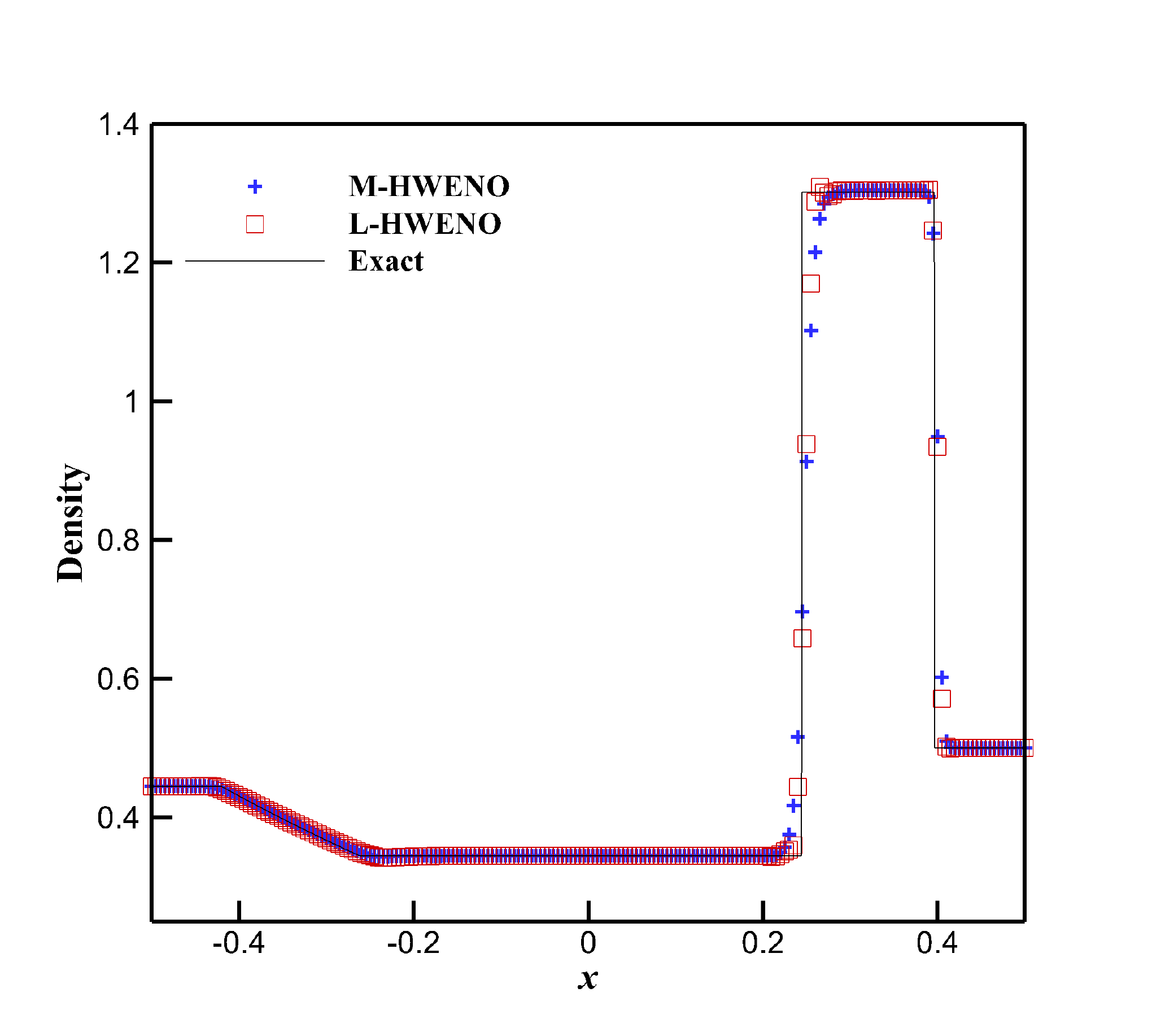}}
\subfigure[close view of (a)]{
\includegraphics[width=0.45\textwidth,trim=20 10 30 30,clip]{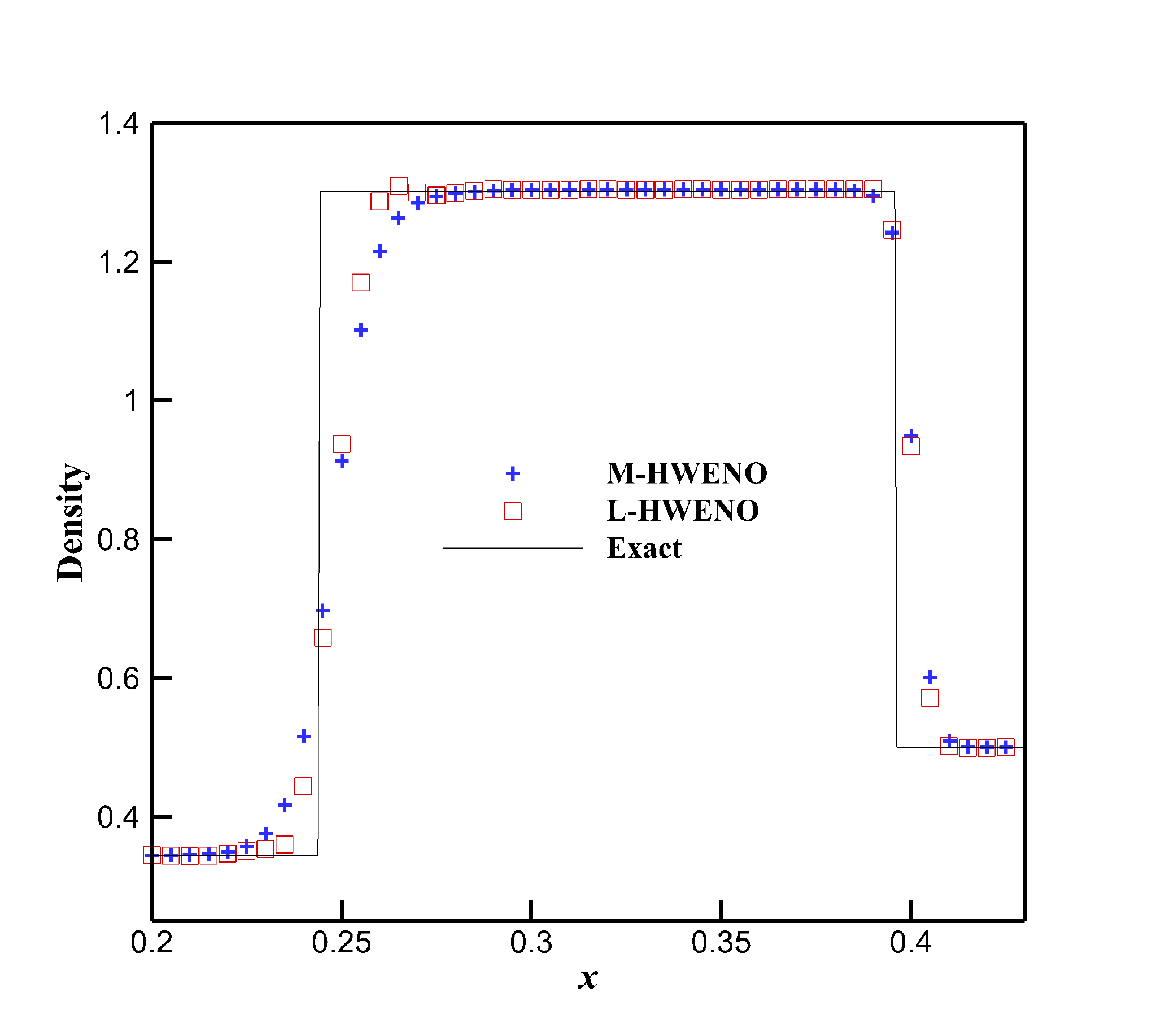}}
 \caption{Example~\ref{lax-1d}. The density $\rho$ at $T=0.16$ obtained by the M-HWENO and L-HWENO schemes.}
\label{laxfig}
\end{figure}

\begin{example}\label{lax-1d}
(Lax problem of the one-dimensional Euler equations.)
\end{example}
In this example, we consider the Lax problem of the one-dimensional Euler equations (\ref{euler1}) with the following initial conditions
\begin{equation*}
\label{lax}
(\rho,\mu,p)= \begin{cases}
(0.445,0.698,3.528), \quad & x<0,\\
(0.5,0,0.571),\quad &  x>0.
\end{cases}
\end{equation*}
The final time is $T=0.16$. The density $\rho$ obtained by the M-HWENO and L-HWENO schemes is presented in Fig. \ref{laxfig}.
We can find that the result obtained by the L-HWENO scheme has a slight higher resolution than that by M-HWENO scheme.
It is worth pointing out that the modification for the derivative of the solution is significant to control oscillations in the M-HWENO scheme \cite{ZhaoZhuang-2020JSC-FD}, while the limiter also plays the same role in the L-HWENO scheme. Similarly, lacking the limiter also generates obvious oscillations, seen in \cite{ZhaoZhuang-2020JSC-FD} for details.

\begin{figure}[H]
\subfigure[density $\rho$: $N_x=400$]{
\includegraphics[width=0.45\textwidth,trim=20 10 30 30,clip]{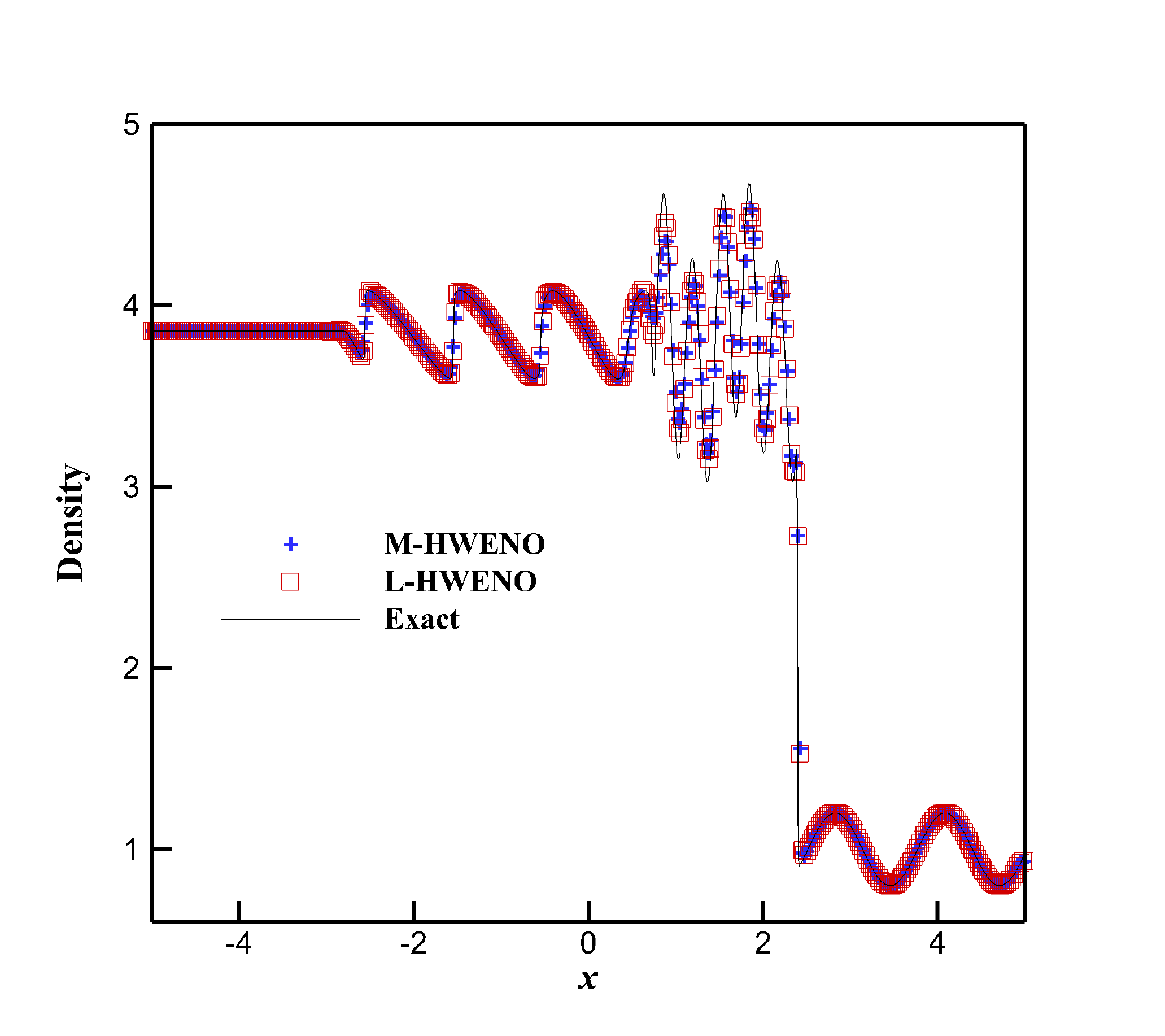}}
\subfigure[close view of (a)]{
\includegraphics[width=0.45\textwidth,trim=20 10 30 30,clip]{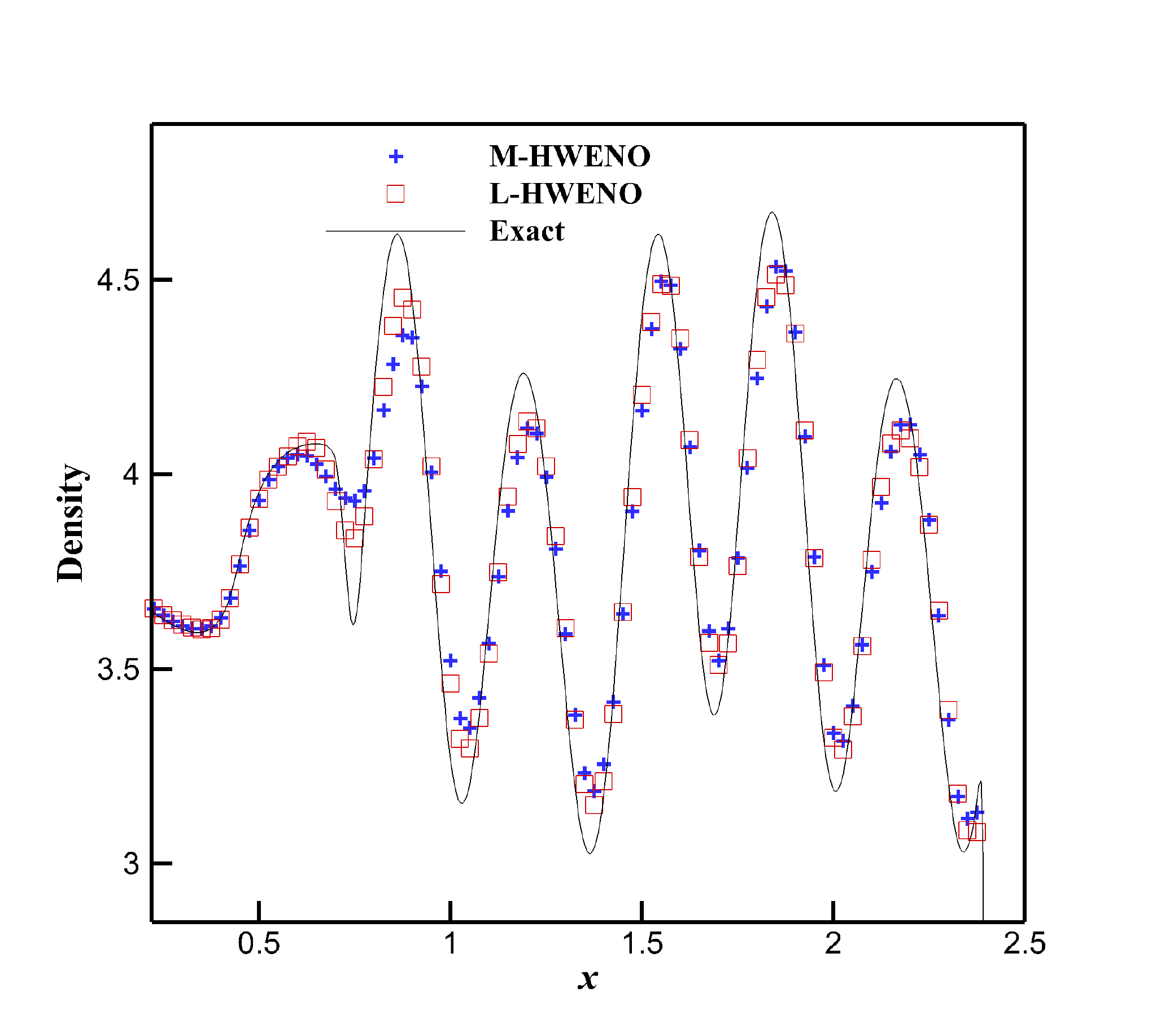}}
\caption{Example \ref{Shu-Osher-1d}. The density $\rho$ at $T=1.8$ obtained by the M-HWENO and L-HWENO schemes.}
\label{sin}
\end{figure}

\begin{example}\label{Shu-Osher-1d}
(Shu-Osher problem of the one-dimensional Euler equations.)
\end{example}
In this example, we consider the following Shu-Osher problem of one-dimensional Euler equations \eqref{euler1}, and the initial conditions are
\begin{equation*}
\label{ShuOsher}
(\rho,\mu,p)= \begin{cases}
(3.857143,~2.629369,~10.333333),\quad & x <-4,\\
(1 + 0.2\sin(5x), ~0, ~1),\quad & x \geq -4.
\end{cases}
\end{equation*}
The final time is $T=1.8$. The solution of this problem has a moving Mach 3 shock interacting with sine waves in density \cite{s3}, and contains both shock waves and complex smooth region structures. The density $\rho$ obtained by the M-HWENO and L-HWENO schemes is shown in Fig.~\ref{sin}, which clearly illustrates that the L-HWENO scheme has higher resolution than the M-HWENO scheme.

\begin{example}\label{blast-1d}
(Two blast waves problem of the one-dimensional Euler equations.)
\end{example}
In this example, we consider a problem of the interaction of two blast waves, and the initial conditions are
\begin{equation*}
\label{blastwave}
(\rho,\mu,p)= \begin{cases}
(1,0,10^3),\quad& 0<x<0.1,\\
(1,0,10^{-2}),\quad& 0.1<x<0.9,\\
(1,0,10^2),\quad& 0.9<x<1.
\end{cases}
\end{equation*}
The final time $T=0.038$ and the reflective boundary condition is applied.

\begin{figure}[H]
\centering
\subfigure[density $\rho$: $N_x=800$]{
\includegraphics[width=0.45\textwidth,trim=20 10 30 30,clip]{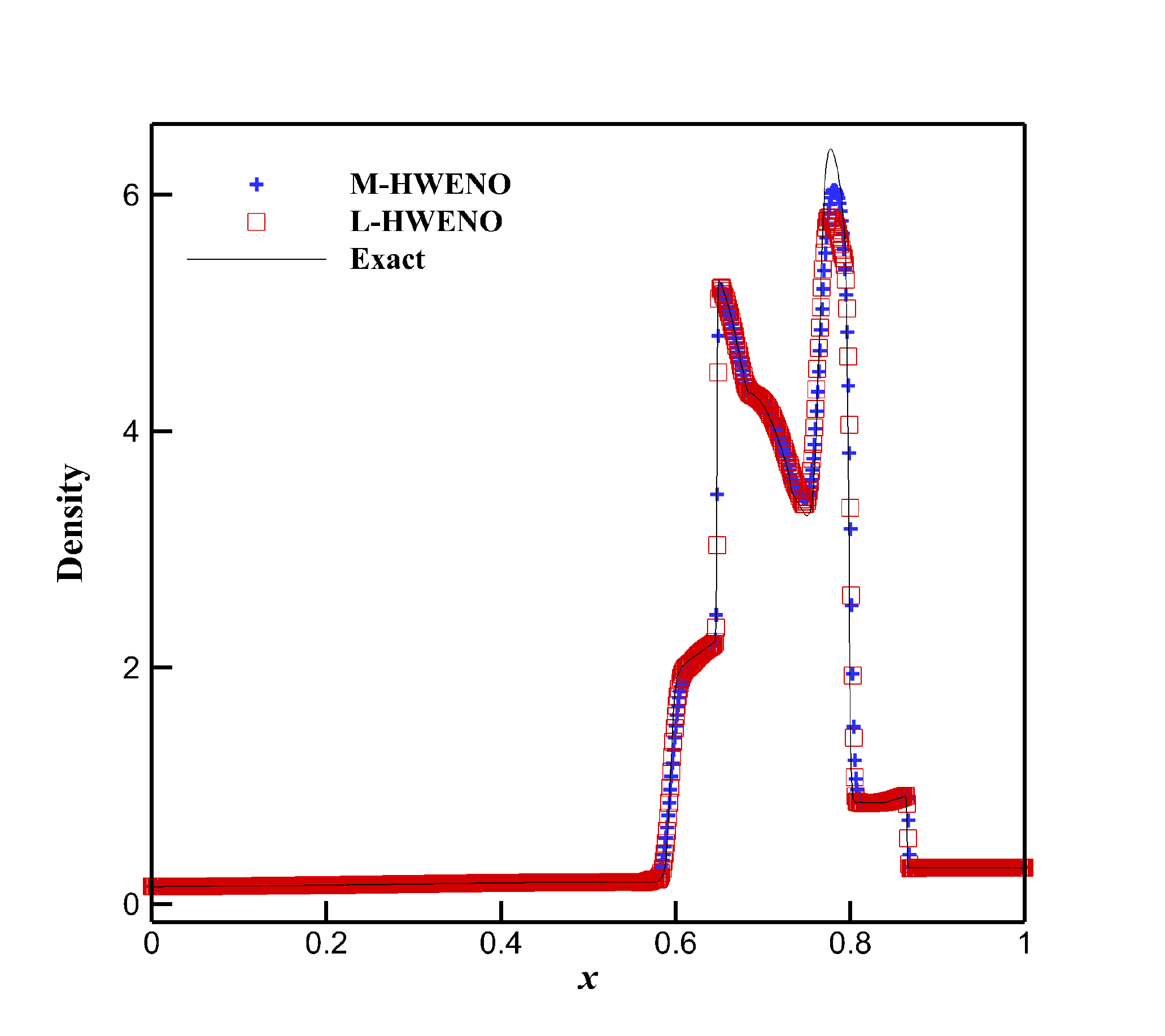}}
\subfigure[close view of (a)]{
\includegraphics[width=0.45\textwidth,trim=20 10 30 30,clip]{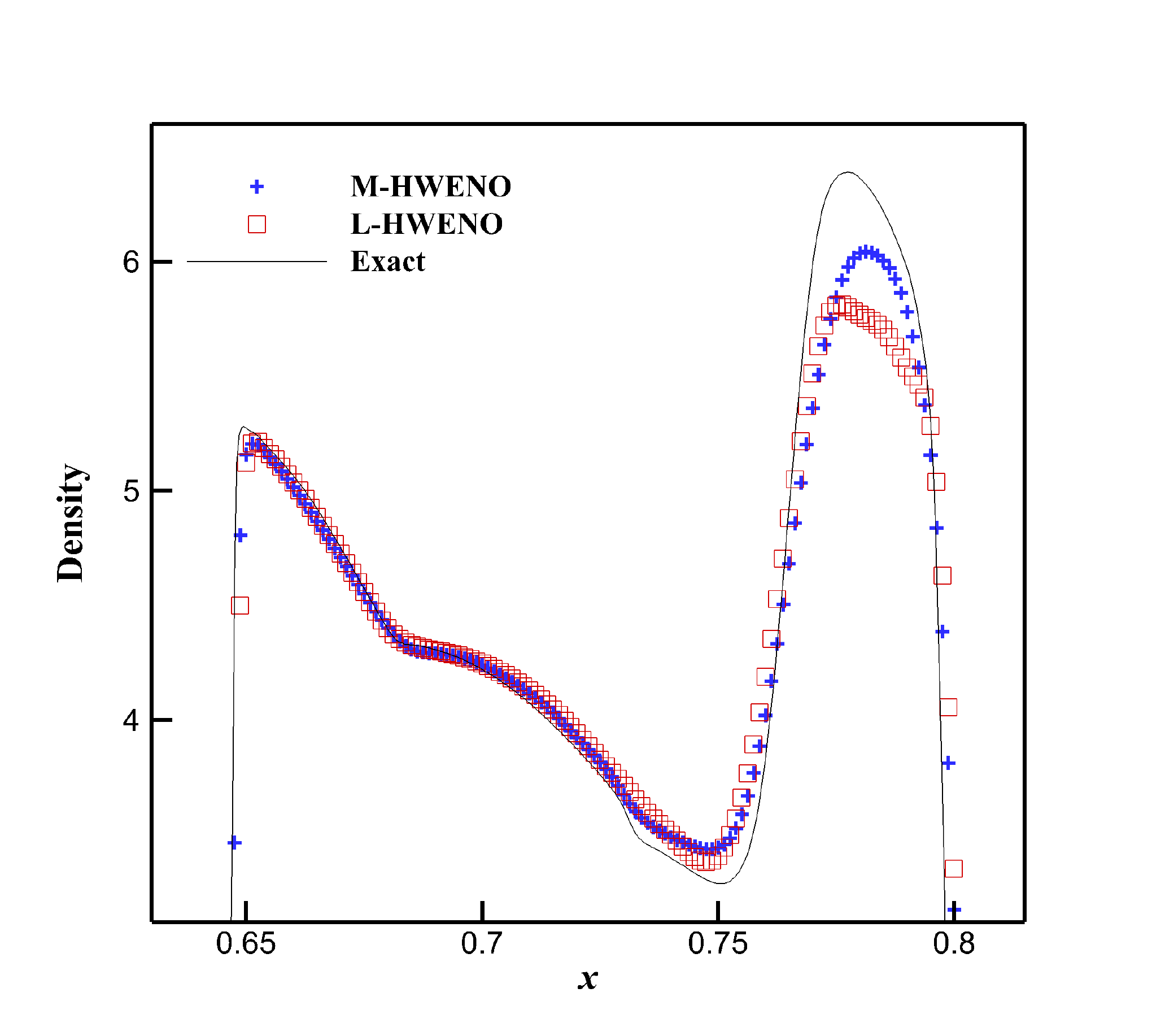}}
\caption{Example \ref{blast-1d}. The density $\rho$ at $T=0.038$ obtained by the M-HWENO and L-HWENO schemes.}
\label{blast1}
\end{figure}

The density $\rho$ obtained by the M-HWENO and L-HWENO schemes at $N_x=800$ is plotted in Fig.~\ref{blast1}.
From the figure, we can know that the resolution of the solution obtained by the L-HWENO scheme  near $x=0.75$ is slightly higher than that obtained by the M-HWENO scheme, but the resolution of the solution obtained by the L-HWENO scheme near $x=0.78$ is slightly lower than that obtained by the M-HWENO scheme.
We think the reason is that two linear polynomials on the small stencils in the limiter have great influence for the derivatives in the L-HWENO scheme. In fact, this phenomenon had been presented in the WENO/HWENO schemes with artificial linear weights \cite{ZQ20,ZQd}.

\begin{figure}[H]
\centering
\subfigure[fixed $G1$ in fluxes reconstruction]{
\includegraphics[width=0.45\textwidth,trim=20 10 30 30,clip]{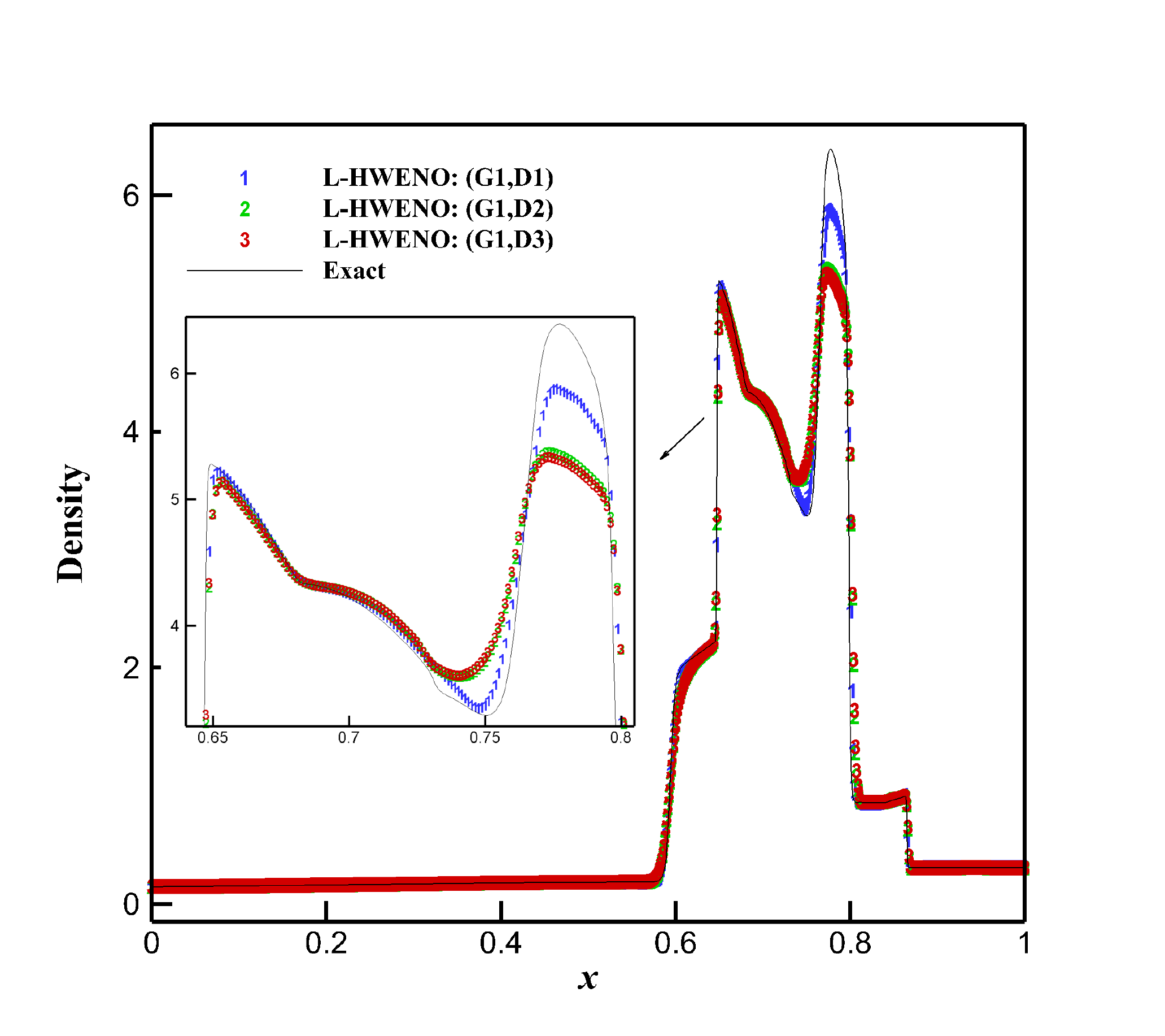}}
\subfigure[fixed $D1$ in limiter]{
\includegraphics[width=0.45\textwidth,trim=20 10 30 30,clip]{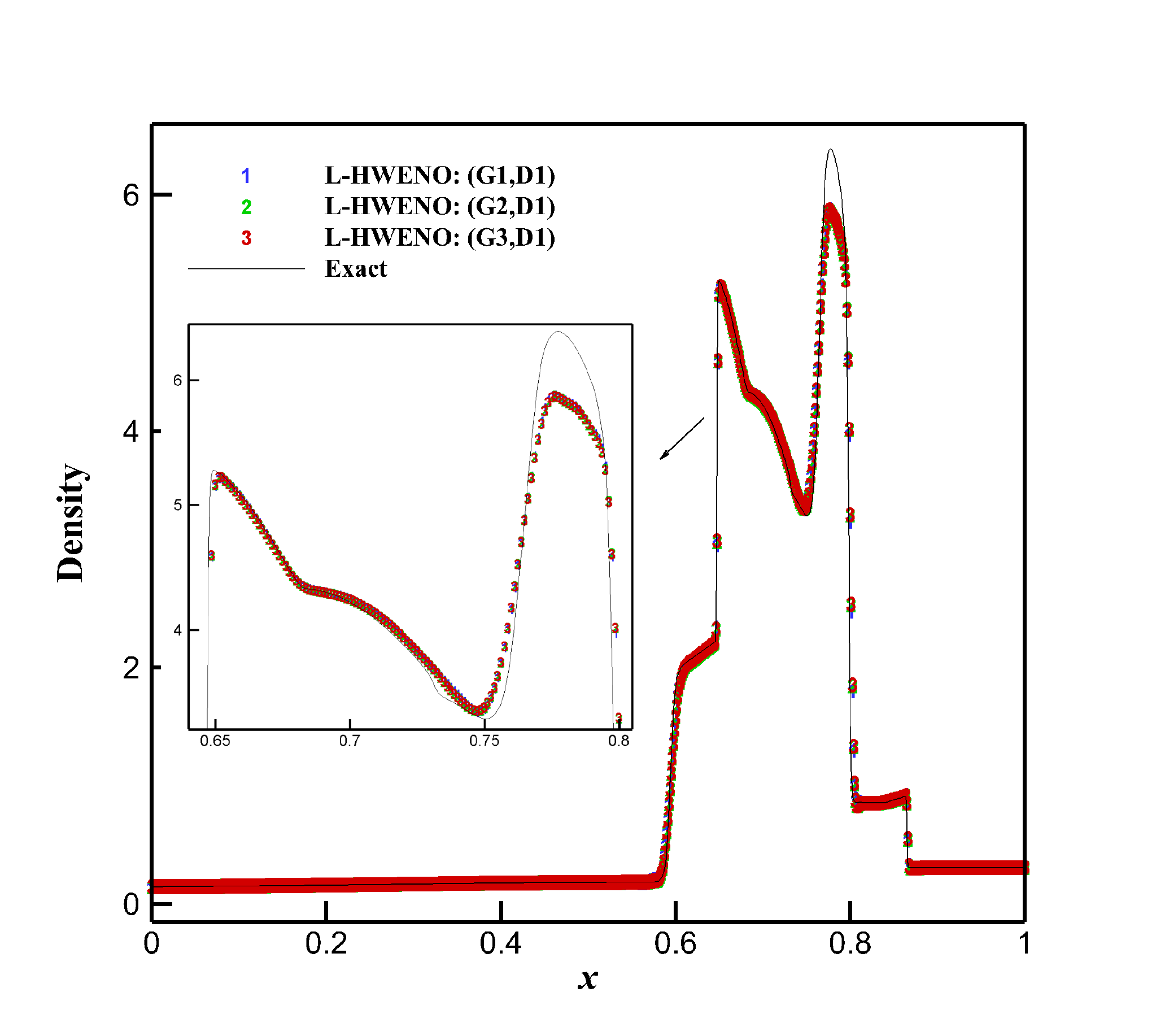}}
\caption{Example \ref{blast-1d}. The density $\rho$ obtained by the L-HWENO scheme at $N_x=800$ with different linear weights setting in limiter and fluxes reconstruction, respectively.}
\label{blast2}
\end{figure}

To study the linear weights how to affect the performance of the L-WENO scheme, we choose several different sets of linear weights $\{\gamma_0,\gamma_1,\gamma_2\}$ (in reconstruction of fluxes) and $\{d_0,d_1,d_2\}$ (in limiter).
The linear weights in the spatial reconstruction are chosen as:
$G1=\{\gamma_0=0.99,~ \gamma_1=\gamma_2=0.005\}$; $G2= \{\gamma_0=\gamma_1=\gamma_2=1/3\}$ and $G3=\{\gamma_0=0.01,~\gamma_1=r_2=0.495\}$.
Similarly, we use the same sets of linear weights in limiter:
$D1=\{d_0=0.99,~ d_1=d_2=0.005\}$; $D2= \{d_0=d_1=d_2=1/3\}$ and $D3= \{d_0=0.01,~d_1=d_2=0.495\}$.
For comparisons, we test the L-HWENO scheme with the linear weights $(G1,D1)$, $(G1,D2)$, $(G1,D3)$, seen in Fig.~\ref{blast2}(a), and combining $(G1,D1)$, $(G2,D1)$ and $(G3,D1)$, seen in Fig.~\ref{blast2}(b).

From Fig.~\ref{blast2}(a), we can find that if the quartic polynomial has larger linear weight in the limiter, the results of the L-HWENO scheme has higher resolution, but it also may have poorer capacity to control non-physical oscillations.
And from the figure \ref{blast2}(b), we know the results are quite similar. That is to say, the resolution of the solution obtained by the L-HWENO scheme is mainly affected by the linear weights in the limiter, and the readers can adjust the linear weights suitably according the explicit problem.
Actually, the L-HWENO scheme combining $(G1,D1)$, $(G2,D1)$ and $(G3,D1)$ has similar performance with the M-HWENO scheme near the peak.

\begin{example}\label{burgers-2d}
(Accuracy test of the two-dimensional Burgers' equation.)
\end{example}
This example is used to verify the accuracy and efficiency of the proposed L-HWENO scheme for the two-dimensional nonlinear Burgers' equation over $[0,4]\times [0,4]$.
The Burgers' equation in two dimensions is
\begin{equation}\label{2dbugers}
  u_t+(\frac {u^2} 2)_x+(\frac {u^2} 2)_y=0.
\end{equation}
The initial condition is $u(x,y,0)=0.5+\sin(\pi (x+y)/2)$ with periodic boundary conditions. The final time is $T=0.5/\pi$ when the solution is still smooth.

\begin{figure}[H]
\centering
\subfigure[Example \ref{burgers-2d}.]{
\includegraphics[width=0.45\textwidth,trim=20 10 30 30,clip]{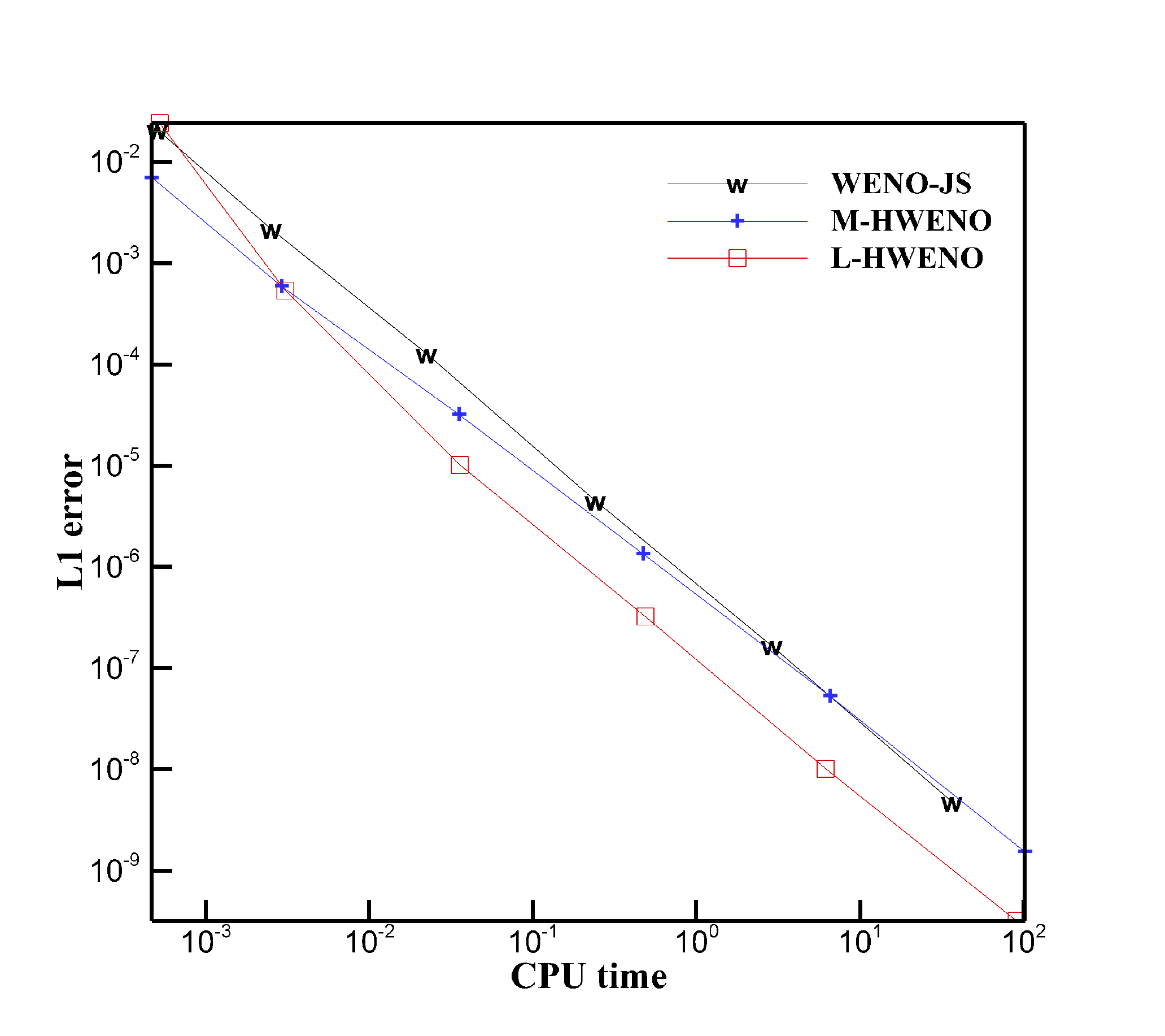}}
\subfigure[Example \ref{Euler-2d}.]{
\includegraphics[width=0.45\textwidth,trim=20 10 30 30,clip]{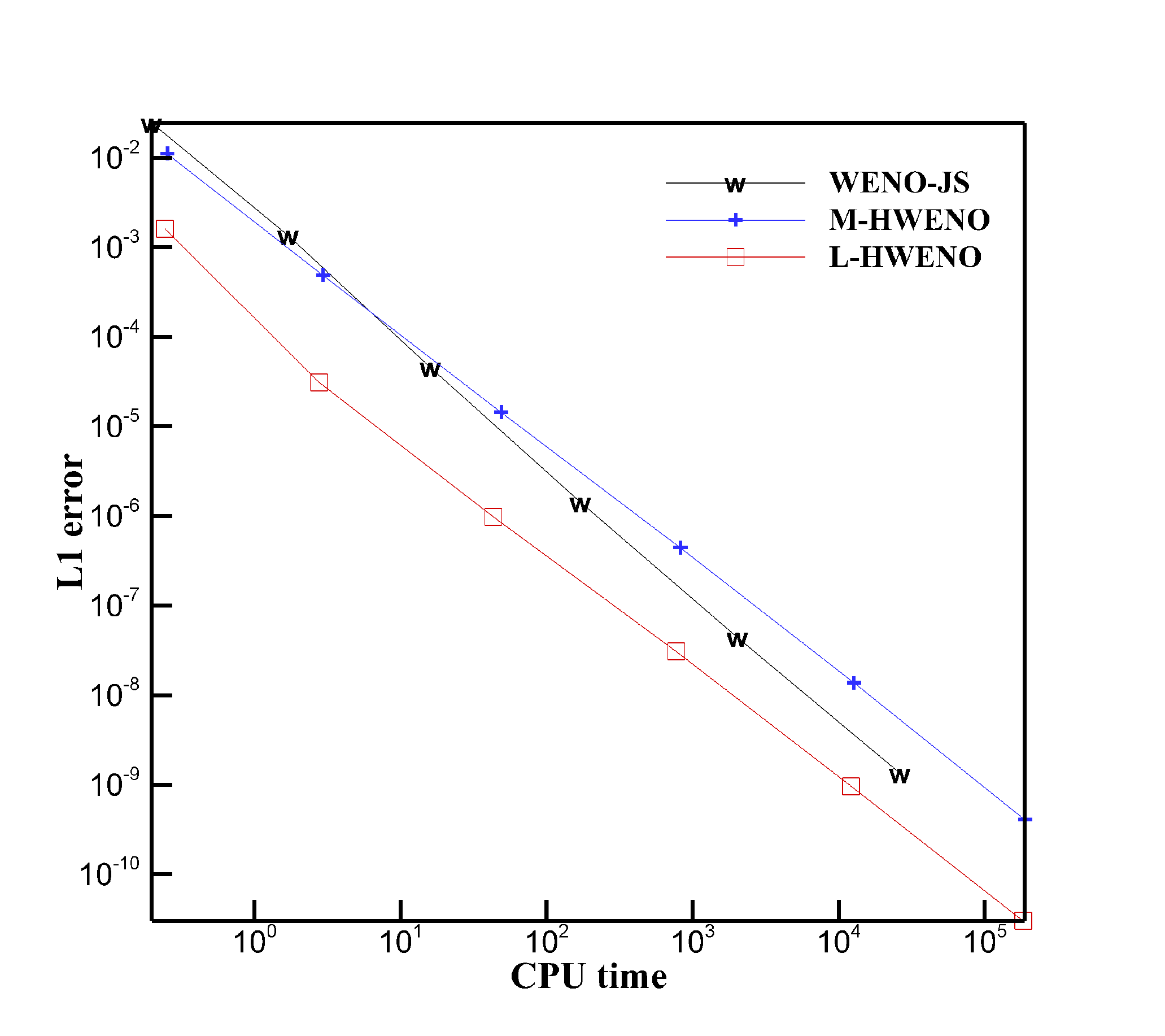}}
 \caption{The error of $L^1$ norm against the CPU time.}
\label{2d-CPU}
\end{figure}

The $L^1$ and $L^\infty$ norm of the error computed by the M-HWENO and L-HWENO schemes are given in Table \ref{tburgers2d}, which
shows the L-HWENO scheme has fifth-order accuracy, and the numerical error of the L-HWENO scheme is smaller than that of M-HWENO scheme. The numerical error against CPU time obtained by the proposed L-HWENO, M-HWENO and WENO-JS schemes is plotted in Fig. \ref{2d-CPU}(a),
which illustrates  the L-HWENO scheme is more efficient than either the M-HWENO scheme or  WENO-JS scheme.

\begin{table}[H]
\centering
\caption{Example~\ref{burgers-2d}. The $L^1$ and $L^\infty$ norm of the error obtained by the M-HWENO and L-HWENO schemes.}
\vspace{5pt}
\begin{tabular} {ccccccccc}
\toprule
 $N_x=N_y$ & \multicolumn{4}{c}{M-HWENO} & \multicolumn{4}{c}{L-HWENO}\\
\cline{2-5} \cline{6-9}
  & $L^1$ error &  order & $L^\infty$ error &  order &
   $L^1$ error &  order & $L^\infty$ error &order\\
      \midrule
$10 $& 6.99E-03 &       & 3.04E-02 &
     & 2.41E-02 &       & 7.30E-02 & \\
$20$ & 5.91E-04 &  3.56 & 3.98E-03 & 2.93
     & 5.33E-04 &  5.50 & 1.90E-03 & 5.26\\
$40$ & 3.22E-05 &  4.20 & 3.82E-04 & 3.38
     & 1.02E-05 &  5.71 & 1.25E-04 & 3.93 \\
$80$ & 1.35E-06 &  4.57 & 1.29E-05 & 4.89
     & 3.23E-07 &  4.98 & 4.28E-06 & 4.86\\
$160$& 5.30E-08 &  4.67 & 5.65E-07 & 4.51
     & 1.02E-08 &  4.99 & 1.34E-07 & 5.00\\
$320$& 1.54E-09 &  5.11 & 2.05E-08 & 4.78
     & 3.19E-10 &  4.99 & 4.25E-09 & 4.97\\
        \bottomrule
\end{tabular}
\label{tburgers2d}
\end{table}

\smallskip

\begin{example}\label{Euler-2d}
(Accuracy test of the two-dimensional Euler equations.)
\end{example}
This example is used to verify the accuracy and efficiency of the proposed L-HWENO scheme for the Euler equations in two dimensions.
The system of Euler equations in two dimensions is
\begin{equation}
\label{euler2}
 \frac{\partial}{\partial t}
 \begin{bmatrix*}[c]
 \rho \\
 \rho \mu \\
 \rho \nu \\
 E
 \end{bmatrix*}
 +
 \frac{\partial}{\partial x}
 \begin{bmatrix*}[c]
\rho \mu \\
\rho \mu^{2}+p \\
\rho \mu \nu \\
\mu(E+p)
\end{bmatrix*}
 +
 \frac{\partial}{\partial y}
 \begin{bmatrix*}[c]
\rho \nu \\
\rho \mu \nu \\
\rho \nu^{2}+p \\
\nu(E+p)
\end{bmatrix*}=0,
\end{equation}
where $\rho$ is the density, $(\mu,\nu)$ is the velocity, $E$ is the total energy and $p$ is the pressure, in which $E = p/(\gamma-1)+\rho(\mu^2+\nu^2)/2 $ with $\gamma=1.4$. The computational domain is $[0,2] \times [0, 2]$. The initial conditions are $
\rho(x,y,0)=1+0.2\sin(\pi(x+y))$, $\mu(x,y,0)=1$, $\nu(x,y,0)=1$, $p(x,y,0)=1$. Periodic boundary conditions are used for all unknown variables.
The exact solution of $\rho$ is $\rho(x,y,t)=1+0.2\sin(\pi(x+y-2t))$.

The final time is $T=2$. The error of the $L^1$ and $L^\infty$ norm computed by the M-HWENO scheme and the proposed L-HWENO scheme is presented in Table \ref{tEluer2d}, which shows the L-HWENO scheme has fifth-order accuracy, and the solution obtained by the L-HWENO scheme is more accurate than that by the M-HWENO scheme.
To show the efficiency of the L-HWENO scheme for this two-dimensional system,
we plot the numerical error against CPU time of the L-HWENO, M-HWENO and WENO-JS schemes in Fig. \ref{2d-CPU}(b). We can find that the error of the L-HWENO scheme is smaller than either the M-HWENO scheme or WENO-JS scheme for a fixed amount of the CPU time.

We also list the $L^1$ and $L^\infty$ norm of the error obtained by the proposed HWENO scheme with and without the limiter in Table \ref{tEluer2d-wo} for this example.
We can clearly see that the proposed HWENO scheme without limiter loses the convergence order at $320\times 320$ points, which violates the common sense because the limiter for the derivatives has the same order accuracy comparing with the reconstruction, therefore, to some extend, the limiter for the derivatives has significant effect to hold stability of the L-HWENO scheme and make its numerical solution be convergent.
The main reason of this phenomenon is the linear approximation of the mixed derivatives in the expression \eqref{mixed}, and this instability also can be solved by splitting the fluxes about the mixed derivatives as other finite different HWENO schemes \cite{LiMRHW1,LiuQ1}, but splitting the fluxes leads to the HWENO schemes \cite{LiMRHW1,LiuQ1} be only the fourth-order accuracy.


\begin{table}[H]
\centering
\caption{Example~\ref{Euler-2d}. The $L^1$ and $L^\infty$ norm of the error computed by the M-HWENO and L-HWENO schemes.}
\vspace{5pt}
\begin{tabular} {ccccccccc}
 \toprule
$N_x=N_y$ & \multicolumn{4}{c}{M-HWENO} & \multicolumn{4}{c}{L-HWENO}\\
\cline{2-5} \cline{6-9}
  & $L^1$ error &  order & $L^\infty$ error &  order
  & $L^1$ error &  order & $L^\infty$ error &  order\\
   \midrule
$     10$ &     1.11E-02 &           &     1.64E-02 &        								&     1.62E-03 &           &     3.57E-03 &          \\
$     20$ &     4.92E-04 &       4.50 &     7.89E-04 &       4.37								&     3.10E-05 &       5.71 &     6.47E-05 &       5.78 \\
$     40$ &     1.44E-05 &       5.09 &     2.69E-05 &       4.87								&     9.87E-07 &       4.97 &     1.64E-06 &       5.30 \\
$     80$ &     4.43E-07 &       5.03 &     8.32E-07 &       5.02								&     3.09E-08 &       5.00 &     4.91E-08 &       5.06 \\
$    160$ &     1.36E-08 &       5.02 &     2.49E-08 &       5.06								&     9.66E-10 &       5.00 &     1.52E-09 &       5.02 \\
$    320$ &     4.08E-10 &       5.06 &     7.00E-10 &       5.16								&     3.02E-11 &       5.00 &     4.75E-11 &       5.00 \\
 \bottomrule
\end{tabular}
\label{tEluer2d}
\end{table}


\begin{table}[H]
\centering
\caption{Example~\ref{Euler-2d}. The $L^1$ and $L^\infty$ norm of the error obtained by the proposed HWENO scheme with and without limiter.}
\vspace{5pt}
\begin{tabular} {ccccccccc}
 \toprule
$N_x=N_y$ & \multicolumn{4}{c}{L-HWENO} & \multicolumn{4}{c}{proposed HWENO without limiter}\\
\cline{2-5} \cline{6-9}
  & $L^1$ error &  order & $L^\infty$ error &  order
  & $L^1$ error &  order & $L^\infty$ error &  order\\
   \midrule
$     10$ &     1.62E-03 &       &     3.57E-03    &       								&     2.16E-03 &           &     3.33E-03 &           \\
$     20$ &     3.10E-05 &       5.71 &     6.47E-05 &       5.78 								&     9.74E-05 &       4.47 &     1.54E-04 &       4.44 \\
$     40$ &     9.87E-07 &       4.97 &     1.64E-06 &       5.30 								&     3.62E-06 &       4.75 &     5.68E-06 &       4.76 \\
$     80$ &     3.09E-08 &       5.00 &     4.91E-08 &       5.06 								&     1.20E-07 &       4.91 &     1.89E-07 &       4.91 \\
$    160$ &     9.66E-10 &       5.00 &     1.52E-09 &       5.02 								&     3.85E-09 &       4.97 &     6.09E-09 &       4.95 \\
$    320$ &     3.02E-11 &       5.00 &     4.75E-11 &       5.00 								&     2.10E-06 &      -9.09 &     1.42E-05 &     -11.18 \\
 \bottomrule
\end{tabular}
\label{tEluer2d-wo}
\end{table}

\begin{figure}[H]
\subfigure[]{
\includegraphics[width=0.45\textwidth,trim=20 10 30 30,clip]{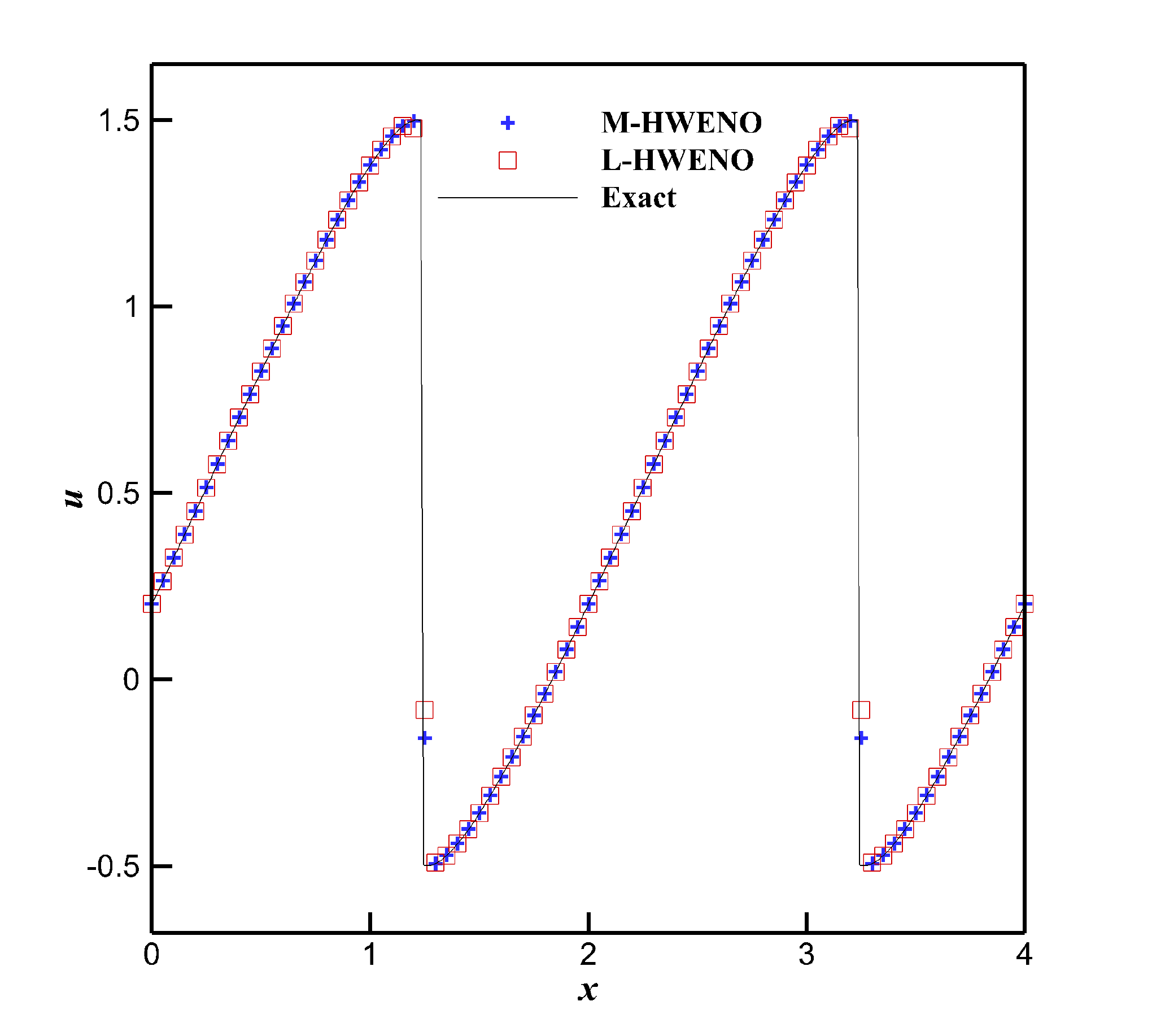}}
\subfigure[]{
\includegraphics[width=0.45\textwidth,trim=20 10 30 30,clip]{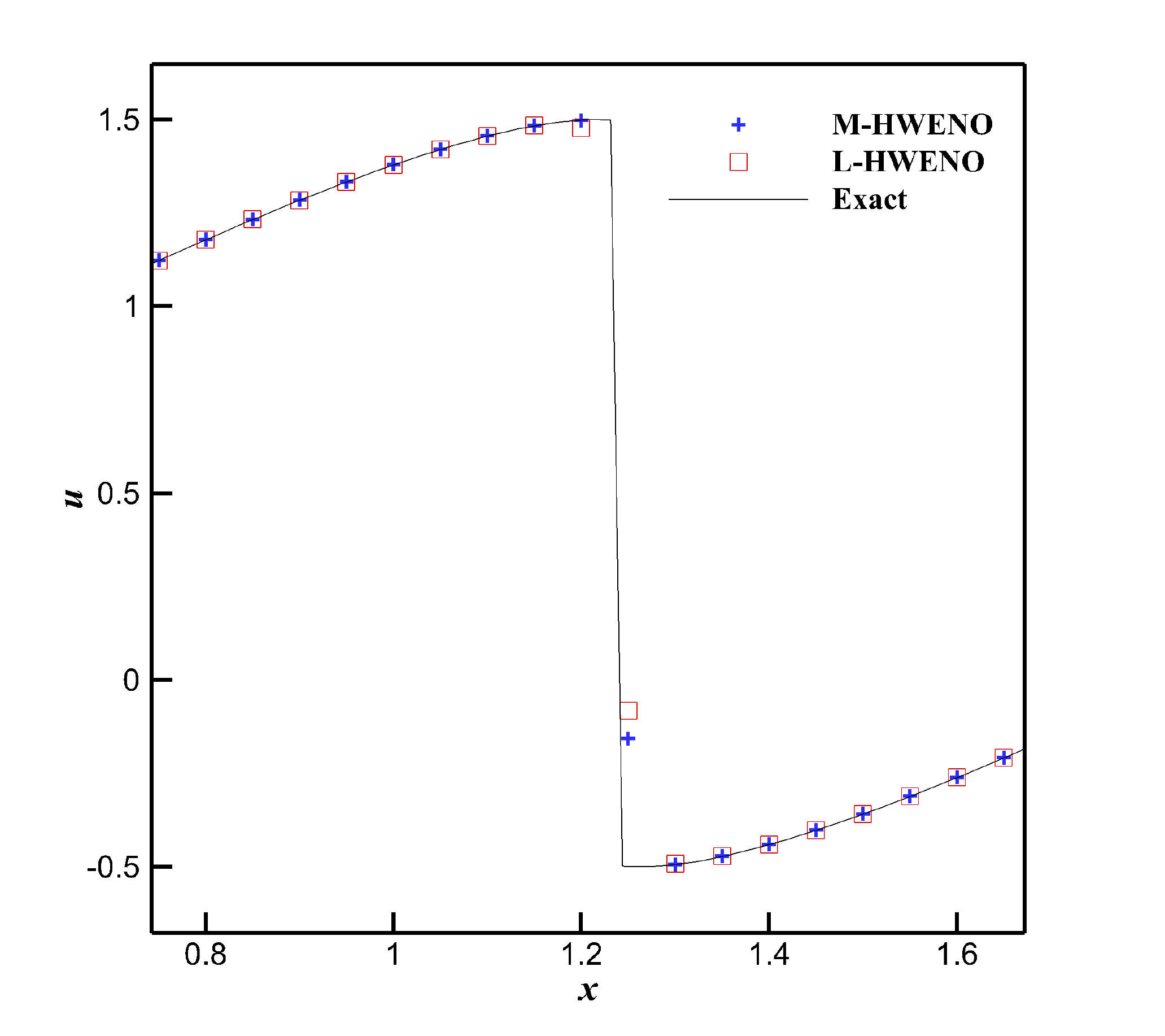}}
 \caption{Example \ref{burgers-2d-shock}. The solution $u$ cut along the line $x=y$ at $T=1.5/\pi$ obtained by the M-HWENO and L-HWENO schemes.}
\label{Fburges2d}
\end{figure}

\begin{example}\label{burgers-2d-shock}
(Shock wave of the two-dimensional Burgers' equation.)
\end{example}
We repeat the two-dimensional Burgers' equation \eqref{2dbugers} given in Example \ref{burgers-2d}, but the final simulation time is $T=1.5/\pi$ when the solution is discontinuous. The numerical solution computed by the M-HWENO and L-HWENO schemes, along with the exact solution is presented in Fig. \ref{Fburges2d}. Again, the two HWENO schemes have similar results with high resolutions.

\begin{example}\label{Double-Mach}
(Double Mach reflection problem of the two-dimensional Euler equation.)
\end{example}
In this example, we test the double Mach reflection problem from \cite{Wooc} modeled by two-dimensional Euler equations (\ref{euler2}) over $[0,4]\times[0,1]$. This example has a reflection wall located at the bottom, starting from $x=\frac{1}{6}$, $y=0$, making a $60^{o}$ angle with the $x$-axis. The exact post-shock condition is imposed from $x=0$ to $x=\frac{1}{6}$ and the rest has the reflection boundary condition for the bottom boundary, and the exact motion of the Mach 10 shock is imposed for the top boundary. Inflow and outflow boundary conditions are used for the left and right boundaries, respectively.

\begin{figure}[H]
\centering
\subfigure[density $\rho$: M-HWENO]{
\includegraphics[width=0.45\textwidth,trim=20 10 40 10,clip]{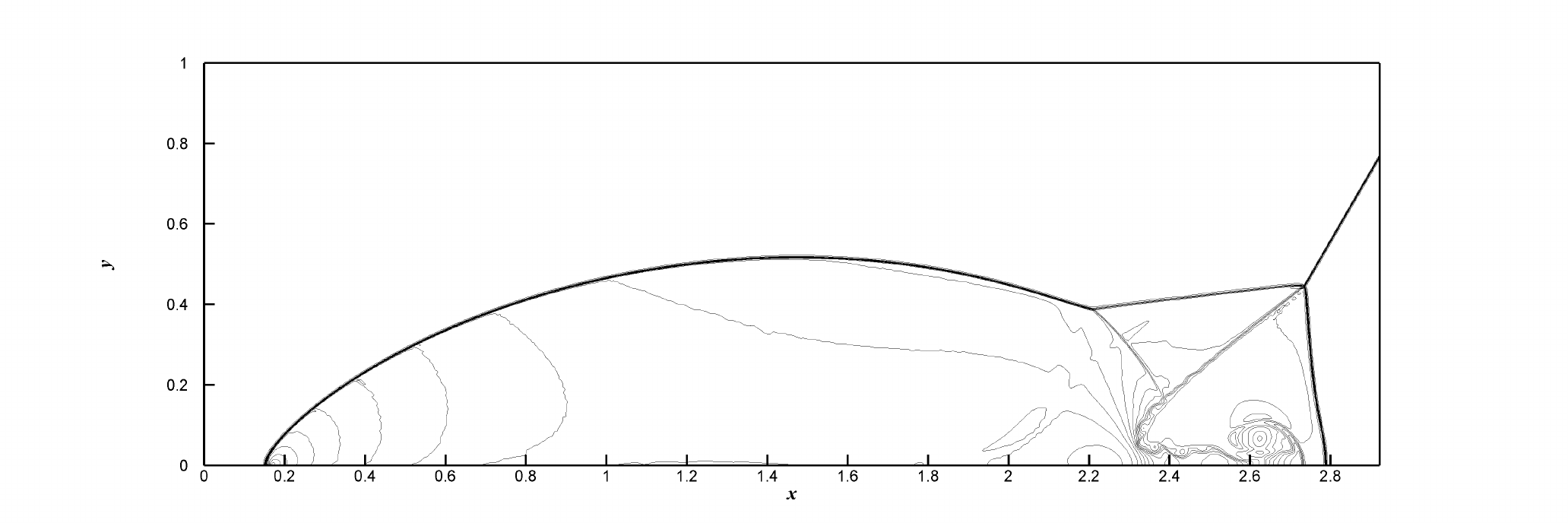} }
\subfigure[density $\rho$: L-HWENO]{
\includegraphics[width=0.45\textwidth,trim=20 10 40 10,clip]{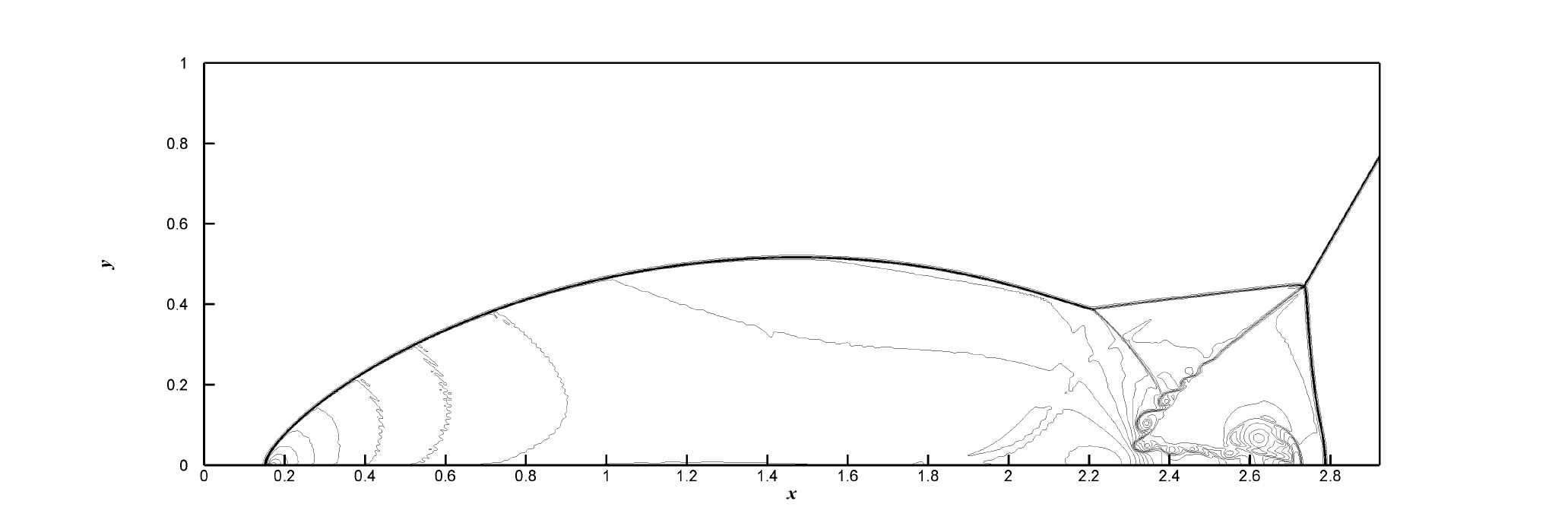} }
\subfigure[close view of (a)]{
\includegraphics[width=0.45\textwidth,trim=0 0 0 0,clip]{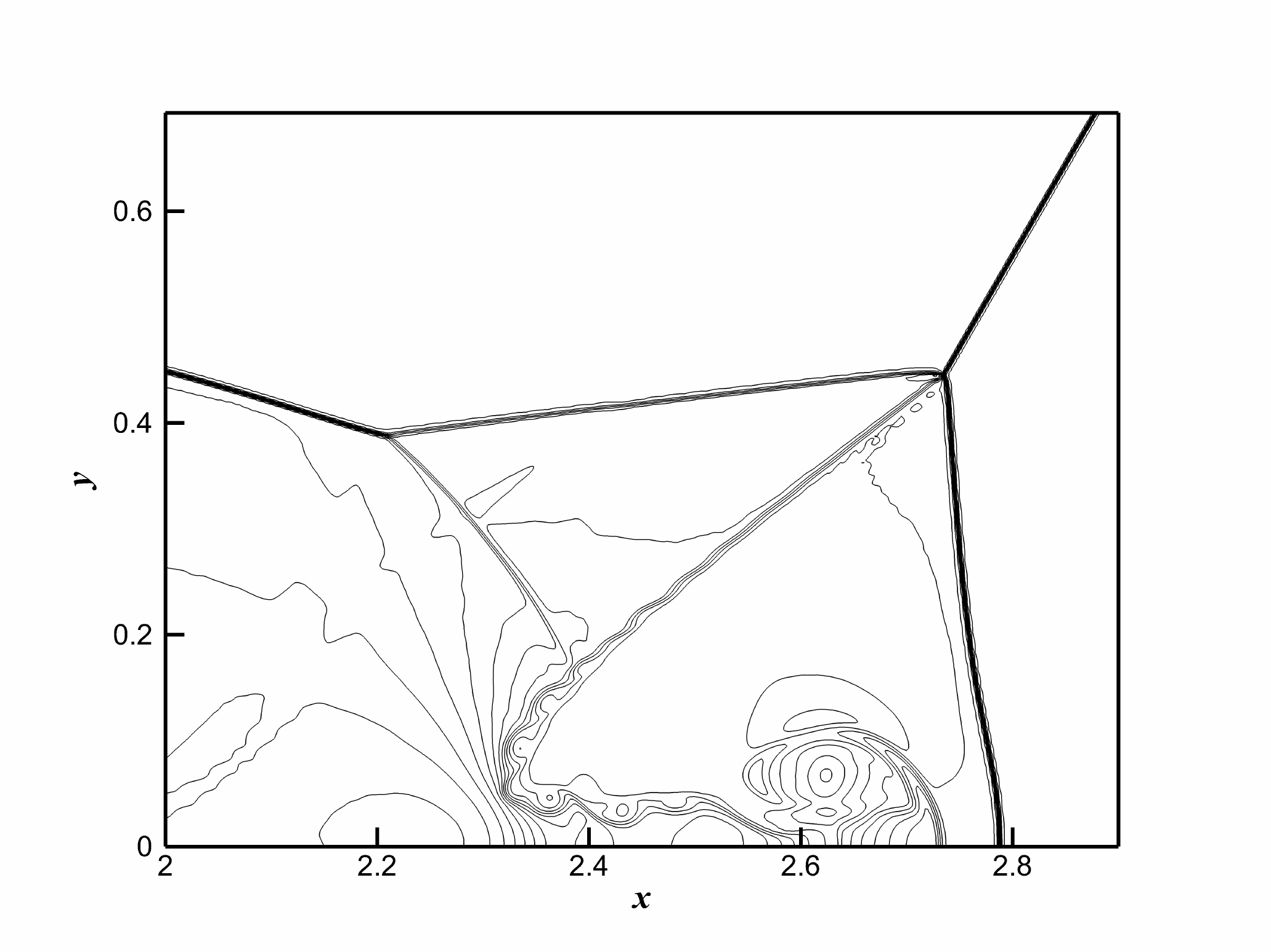}}
 \subfigure[close view of (b)]{
\includegraphics[width=0.45\textwidth,trim=0 0 0 0,clip]{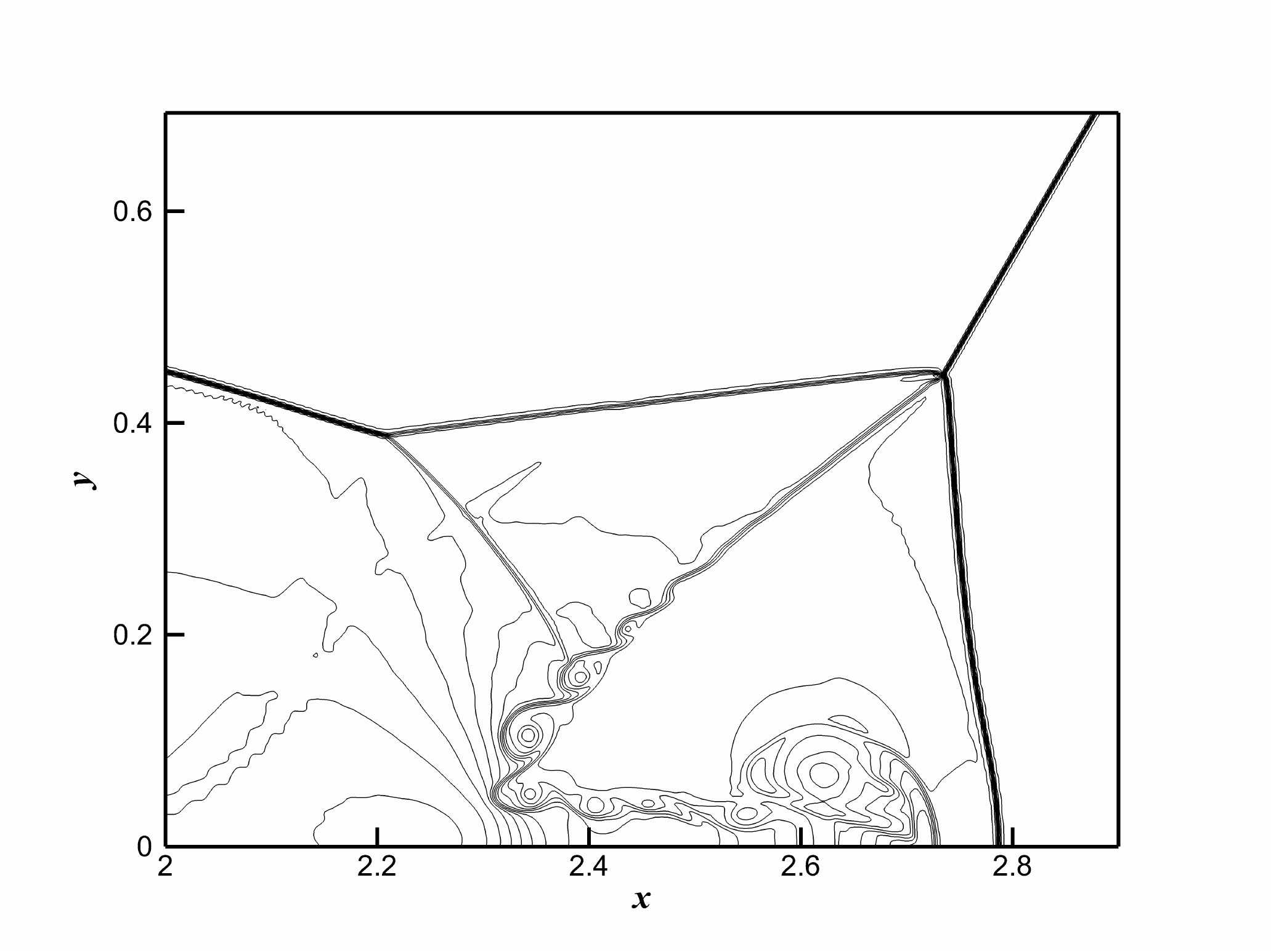} }
 \caption{Example \ref{Double-Mach}. The density $\rho$ at $T=0.2$ obtained by the M-HWENO and L-HWENO schemes.}
\label{smhfig}
\end{figure}

The final time is $T=0.2$. In Fig.~\ref{smhfig}, we show the numerical results computed by the M-HWENO and  L-HWENO schemes in the region $[0,3]\times[0,1]$ and the blow-up region around the double Mach stems. It is observed that the L-HWENO scheme has higher resolution than the M-HWENO scheme, and the L-HWENO scheme captures more complicated structures.
\smallskip

\begin{example}\label{Forward-step}
(Forward step problem of the two-dimensional Euler equation.)
\end{example}
In this example, we test a forward step \cite{Wooc} modeled of the two-dimensional Euler equations (\ref{euler2}).
There is a wind tunnel with a initial right-going Mach 3 flow, and it has the width of 1 unit and the length of 3 units. The location of the step corner is $(x,y)=(0.6,0.2)$.
Reflective boundary conditions are used along the wall of the tunnel. Inflow and outflow boundary conditions are used at the entrance and the exit, respectively. The corner of the step is a singular point and we treat it as in \cite{Wooc}.

\begin{figure}[H]
\centering
\subfigure[density $\rho$: M-HWENO]{
\includegraphics[width=0.45\textwidth,trim=10 10 10 0,clip]{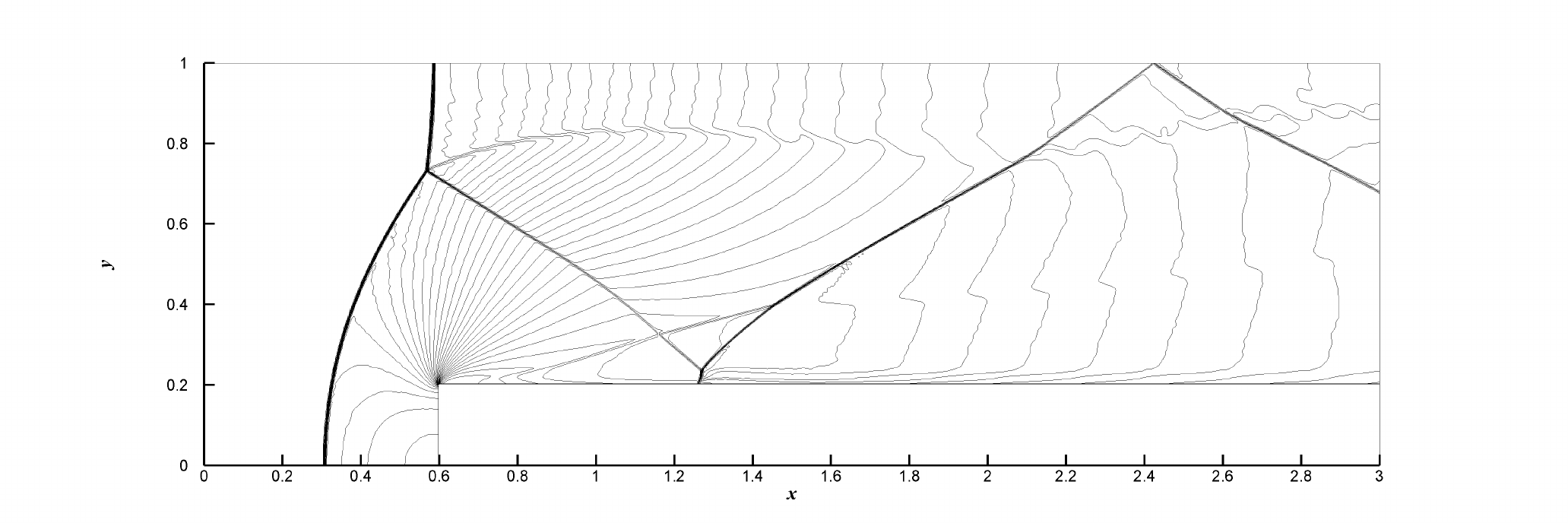} }
\subfigure[density $\rho$: L-HWENO]{
\includegraphics[width=0.45\textwidth,trim=10 10 10 0,clip]{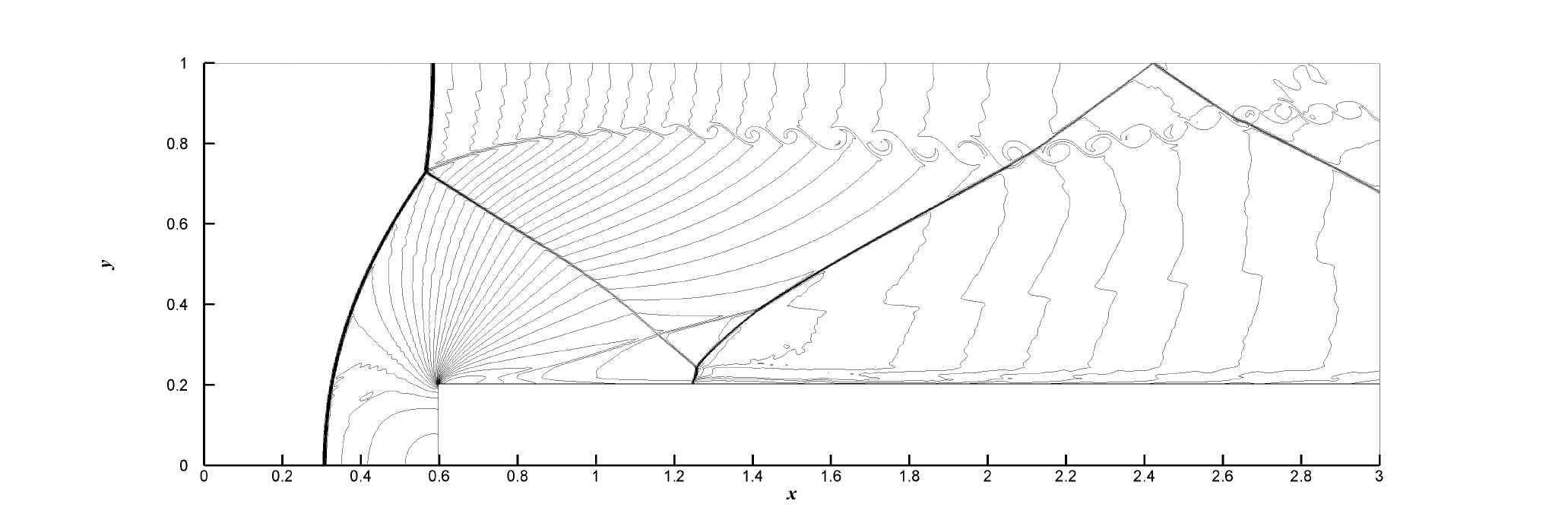}  }
 \caption{ Example \ref{Forward-step}. The density $\rho$ at $T=4$ obtained by the M-HWENO and L-HWENO schemes.}
\label{stepfig}
\end{figure}

We compute the time up to $T=4$. The numerical results of the M-HWENO and L-HWENO schemes at $960\times320$ grid points are shown in Fig. \ref{stepfig}. We can observe that the L-HWENO scheme has higher resolution than the M-HWENO scheme.

\section{Conclusions}
\label{sec:conclusions}

In this paper, a simple fifth-order finite difference Hermite weighted essentially non-oscillatory (HWENO) scheme combined with limiter (called as the L-HWENO scheme) is constructed for one- and two- dimensional hyperbolic conservation laws.
The fluxes in the governing equation are approximated by the nonlinear HWENO reconstruction which is the combination of a quintic polynomial with two quadratic polynomials, where the linear weights can be artificial positive numbers as long as their sum equals one. And other fluxes are approximated by high-degree polynomials directly, which leads to the result that the reconstruction of the fluxes for derivative equations is linear. For the purpose of controlling spurious oscillations, an HWENO limiter is applied to modify the derivatives as
the modified HWENO (M-HWENO) scheme \cite{ZhaoZhuang-2020JSC-FD}. Instead of using the modified derivatives both in fluxes reconstruction and time discretization as in \cite{ZhaoZhuang-2020JSC-FD}, we only apply the modified derivatives in time discretization while remaining the original derivatives in fluxes reconstruction.
Comparing with the M-HWENO scheme \cite{ZhaoZhuang-2020JSC-FD}, the proposed L-HWENO scheme is simpler, more accurate, efficient, and higher resolution.

The spatial reconstruction and the limiter for the derivatives both use a high-degree polynomial combined with two lower-degree polynomials convexly, where the corresponding linear weights can be any positive numbers (their sum is 1).  It is easy to implement and has the ability to capture complicated structures. In the implementation, the limiter in the proposed L-HWENO scheme plays an important role to improve stability and keep  high resolution, where lacking this procedure would lead to instability in two dimensions even for a linear problem (cf. Example \ref{Euler-2d}). Meanwhile, different linear weights in the limiter will impact the resolution near discontinuities, while the linear weights in the spatial reconstruction have a slight effect (cf. Example \ref{blast-1d}).

Various benchmark numerical examples have been tested to demonstrate the accuracy and efficiency of the L-HWENO scheme. The results show that the L-HWENO scheme has fifth-order accuracy, and the solution of the L-HWENO scheme is more accurate than that of the M-HWENO scheme. Meanwhile, the L-HWENO scheme is more efficient than either the M-HWENO scheme \cite{ZhaoZhuang-2020JSC-FD} or WENO-JS scheme \cite{js}.
Note that the efficiency of the M-HWENO scheme is slightly less than that of the WENO-JS scheme. In addition, we would mention that the L-HWENO and M-HWENO schemes both use a compact three-point reconstructed stencil while a five-point stencil is need in \cite{js} even though they all have the fifth-order accuracy.

\end{document}